\documentclass[12pt]{article}

\usepackage{fancyheadings}
\usepackage{color}
\usepackage{latexsym}
\usepackage{setspace}
\usepackage{amsmath}
\usepackage{amssymb}
\usepackage{latexsym}
\usepackage{amsmath}
\usepackage{amssymb}
\usepackage{natbib}
\usepackage{subfigure}
\usepackage{multirow}
\usepackage{colortbl}
\setcitestyle{authoryear,open={(},close={)}}
\usepackage[table,xcdraw]{xcolor}
\usepackage{multirow}
\setcitestyle{numbers}
\setcitestyle{square}
\usepackage{hyperref}
\hypersetup{colorlinks,citecolor=blue}
\usepackage{graphics}
\usepackage{graphicx}
\makeatletter
\def\ps@pprintTitle{%
  \let\@oddhead\@empty
  \let\@evenhead\@empty
  \let\@oddfoot\@empty
  \let\@evenfoot\@oddfoot
}
\makeatother
\makeatletter \oddsidemargin  -.5cm \evensidemargin -.5cm \textwidth
18cm \topmargin -.60in \textheight 24cm

\newtheorem{theorem}{Theorem}[section]

\newtheorem{lemma}[theorem]{Lemma}

\newtheorem{definition}{Definition}[section]
\newtheorem{corollary}[theorem]{Corollary}
\newtheorem{remark}{Remark}[section]
\newtheorem{example}{Example}[section]

\numberwithin{equation}{section}
\numberwithin{table}{section}
\numberwithin{equation}{section}
\begin{document}

\title{Estimating location parameters of two exponential distributions with ordered scale parameters}
\author{\small $^a$Lakshmi Kanta Patra\footnote
	{\baselineskip=10pt~patralakshmi@gmail.com, lkpatra@iitbhilai.ac.in}, $^{b}$Constantinos Petropoulos, $^a$Shrajal Bajpai and $^c$ Naresh Garg \\
    \small $^{a}$Department of Mathematics, Indian Institute of Technology Bhilai\\
    \small $^{b}$Department of Mathematics, University of Patras, Rio, Patras, Greece 26500\\
\small $^{c}$Statistics and Mathematics unit, Indian Statistical Institute Delhi}
\spacing{1.1}
\date{}
\maketitle
\vspace{-0.5cm}
\begin{abstract}
In the usual statistical inference problem, we estimate an unknown parameter of a statistical model using the information in the random sample. A priori information about the parameter is also known in several real-life situations. One such information is order restriction between the parameters. This prior formation improves the estimation quality. In this paper, we deal with the component-wise estimation of location parameters of two exponential distributions studied with ordered scale parameters under a bowl-shaped affine invariant loss function and generalized Pitman closeness criterion. We have shown that several benchmark estimators, such as maximum likelihood estimators (MLE), uniformly minimum variance unbiased estimators (UMVUE), and best affine equivariant estimators (BAEE), are inadmissible. We have given sufficient conditions under which the dominating estimators are derived. Under the generalized Pitman closeness criterion, a Stein-type improved estimator is proposed. As an application, we have considered special sampling schemes such as type-II censoring, progressive type-II censoring, and record values. Finally, we perform a simulation study to compare the risk performance of the improved estimators.
\noindent {\color{blue} \textit{MSC2020:} 62F10, 62F30, 62C99}
\end{abstract}

\section{Introduction}
Statistical inference problems under order-restricted parameters are extensively studied in the literature because of their application in several applied fields such as clinical trials, bio-assays, genetics, reliability, life testing studies, etc.
A detailed review of the earlier work related to statistical inference under order restrictions can be found in \cite{barlow1972statistical}, \cite{robertson1988}, and \cite{silvapulle2004constrained}. 

The two-parameter exponential distribution is extensively used in life testing and reliability engineering. The problem of estimating location parameters of exponential distributions with ordered scale parameters is encountered in many real-life applications. In reliability theory, the exponential distribution is used to model the lifetime of an experimental unit. The location parameter of an exponential distribution represents the component's guarantee time, and the reciprocal of the scale parameter represents the instantaneous failure rate or hazard rate. So, order relation between the scale parameters occurs naturally in life testing experiments. Suppose a system has two components, C1 and C2. The lifetime of C1 follow $Exp(\mu_1,\sigma_1)$ and lifetime of C2 follow $Exp(\mu_2,\sigma_2)$. Then, it is natural to consider the order relation between the scale parameters $\sigma_1$ and $\sigma_2$. In this article, we will study component-wise estimation of location parameters $\mu_1$ and $\mu_2$ under the restriction $\sigma_1 \le \sigma_2$ under a bowl-shaped affine invariant loss function. In the existing literature, several researchers studied the problem of estimating the ordered parameters. Some related work in this direction as follows  
\cite{vijayasree1991simultaneous}, \cite{pal1992order},   \cite{kumar1988simultaneous}, \cite{shi1998maximum} and reference their in.
\\
Let $(X_{i1},\ldots, X_{in_{i}})$ be a random sample drawn from the $i$-th population $\mathcal{P}_{i},$ $i=1,2$. The probability density function of the $i$-th population $\mathcal{P}_{i}$ is given by
\begin{eqnarray}
	f_{i}(x;\mu_{i},\sigma_{i})=\left\{\begin{array}{ll}
		\frac{1}{\sigma_{i}}~\exp{\Big(-\frac{x-\mu_{i}}{\sigma_{i}}\Big)},
		& \textrm{if $x>\mu_{i},$}\\
		0,& \textrm{otherwise,}
	\end{array} \right.
\end{eqnarray}
where $-\infty<\mu_{i}<\infty$ and $\sigma_{i}>0.$ We assume that
$\sigma_{i}$'s and $\mu_{i}$'s, $i=1,2$ are unknown. Further, we
assume an additional condition on the scale parameters as
$\sigma_{1}\leq\sigma_{2}$. Let us denote
\begin{eqnarray*}
	X_{i(1)}=\min_{1\leq j\leq
		n_{i}}X_{ij},~~~~T_{i}=\sum_{j=1}^{n_{i}}\left(X_{ij}-X_{i(1)}\right),
\end{eqnarray*}
where $i=1,2.$ It is known that
$\left(\underline{X_{(1)}},\underline{T}\right)$ is a complete and
sufficient statistic, where
$\underline{X_{(1)}}=(X_{1(1)},X_{2(1)})$ and
$\underline{T}=(T_{1},T_{2})$. Moreover,
$(X_{1(1)},X_{2(1)},T_{1},T_{2})$ are independently distributed with
\begin{equation} \label{model1}
X_{i(1)}\sim E\left(\mu_{i},\frac{\sigma_{i}}{n_{i}}\right)~~\mbox{ and }~~~
T_{i}\sim Gamma(n_{i}-1,\sigma_{i}), i=1,2. 
\end{equation}
Our aim is to estimate $\mu_i$ respect to a bowl shaped affine invariant loss function $L\left(\frac{\delta-\mu_i}{\sigma_i}\right)$. We assume that the loss function satisfies the following conditions:
\begin{itemize}
	\item [(C1)] $L(t)$ is a twice differentiable strictly bowl-shaped function and non-monotone.
	\item [(C2)] Integrals involving $L(t)$ are finite, and differentiation under the integral sign is permissible.
\end{itemize}
Now, we will give a literature review for estimating ordered parameters of exponential distributions. \cite{vijayasree1991simultaneous} have considered the mixed estimator to estimate the ordered means of two exponential distributions. They proved that mixed estimators dominate the usual estimator. \cite{misra1994estimation} discussed the problem of estimating the ordered location parameter of two exponential distributions with scale parameters are known. They have obtained a class of admissible estimators in the class of mixed estimators. They have also proved the inadmissibility of unrestricted MRE under the order restriction.   \cite{vijayasree1995componentwise} considered the $k (\ge2)$ exponential distributions and investigated the component-wise estimation of the parameters under order restriction. For this estimation problem, they have considered location-invariant quadratic loss.  \cite{vijayasree1995componentwise} used the techniques of \cite{brewster1974improving},  given sufficient condition for the inadmissibility of unrestricted MRE.   Improved estimation of scale parameters of two gamma distributions under ordered restriction has been studied by  \cite{misra2002smooth}. They obtained smooth improved estimators using the techniques of \cite{brewster1974improving}. \cite{chang2008estimation}  discussed the problem of estimating the linear function of ordered scale parameters of two gamma populations under the entropy loss function. \cite{jana2015estimation} have investigated the ordered scale parameters of two exponential distributions with a common location parameter. They have obtained  UMVUE and restricted MLE. \cite{jana2015estimation} also derived improved estimators in classes of scale and affine equivariant estimators. Component-wise estimation of ordered scale parameter two populations of the multivariate Lomax distribution investigated by \cite{petropoulos2017estimation}. The author proved the inadmissibility of best equivariant estimators and also constructed estimators that dominate usual estimators.   \cite{patra2018estimating} have considered estimating the reciprocal of common scale parameters with unknown and ordered location parameters. They showed that the best equivariant estimator is inadmissible. Recently \cite{patra2021minimax} investigated the problem of estimating common variance and precision of two normal distributions with ordered means. They have proposed several improved estimators. \cite{patra2021componentwise} discussed the component-wise estimation of ordered scale parameters of two experiential distributions under a bowl-shaped scale invariant loss function. Recently \cite{garg2022estimation} have studied component-wise estimation of ordered location parameters of a bivariate location family and bivariate under order restriction with respect to the general loss function. They have given sufficient conditions under which improved estimators are derived. 

Similar to the estimation with respect to the loss function, various researchers have derived an improved estimator based on the pitman closeness criteria. \cite{nayak1990estimation} have discussed the estimation of location and scale parameters through the application of the generalized Pitman nearness criterion. Extensive research has been conducted using the Pitman closeness approach. For in-depth exploration of the early work in this field, one can consult \cite{robert1993pitman}, \cite{das1991comparison}. However, there are only a limited number of studies on the use of the PN criterion in estimating order-restricted parameters. In recent work, \cite{ahmadi2010pitman} formulated expressions to determine the Pitman closeness of current records to a common population parameter. \cite{balakrishnan2012pitman} have evaluated the performance of the best linear unbiased and best linear invariant predictors for censored order statistics in the one-sample scenario and order statistics from a future sample in the two-sample scenario, using the pitman closeness criterion by assuming a type-II right-censored sample originating from an exponential distribution. \cite{balakrishnan2013pitman} have done comparisons of Pitman closeness for type-I censored samples from an exponential distribution. \cite{chang2015estimation} have investigated the problem of estimating two ordered normal means  using the modified Pitman nearness criterion, both when the order restriction on variances is present and when it is not. \cite{mirmostafaee2015pitman} obtained the optimal equivariant predictor and conducted a comparison with the best invariant predictor and the best unbiased predictor. \cite{garg2023unified} performed component-wise estimation of location and scale parameters of a general bivariate location and scale probability distribution subject to order restriction  under the Generalized Pitman Nearness (GPN) standard, along with a general loss function. \cite{garg2023improved} investigated the problem of estimating the larger location parameter within two general location families, employing criteria that involve minimizing the risk function and achieving Pitman closeness, all within the context of a general bowl-shaped loss function.

In the present work, we have extended the work of \cite{vijayasree1995componentwise} for a general bowl-shaped affine invariant loss function with $k=2$ and for generalized Pitman closeness criterion. We are interested in proving the inadmissibility of well-known estimators such as MLE, BAEE and UMVUE. To do so, we will obtain various  dominating estimators. We have also proposed generalized Bayes estimators.

The rest of the paper is organized as follows. In Section \ref{sec2s}, we have discussed the estimation of $\mu_1$. We propose estimators with a smaller risk than the BAEE, MLE, and UMVUE of $\mu_1$ under a bowl-shaped affine invariant loss function. The improved estimators are derived for squared error and linex loss function as an application. Also, we have proposed a class of improved estimators using the IERD approach of \cite{kubokawa1994unified}. We have seen that the boundary estimator of this class is a generalized Bayes estimator. Estimation of $\mu_2$ will be studied in Section \ref{sec3}. We derive results similar to Section \ref{sec3}  for estimating $\mu_2$.  We have discussed the improved estimation of $\mu_1$ and $\mu_2$ under generalized Pitman closeness criterion in Section \ref{sec4}. In Section \ref{sec5} we have considered improved estimation based on special sampling schemes such as type-II censoring, progressive type-II censoring and record values.  A simulation study has been carried out in Section \ref{sec6} to compare the risk performance of the improved estimators.    
\section{Improved estimation of $\mu_1$}\label{sec2s}
Here we will discuss estimation of the location parameter  $\mu_{1}$ with $\sigma_{1}\leq\sigma_{2}$. At first we invoke the principle of invariance. For this purpose  consider the affine group of transformations
\begin{eqnarray}\label{eqtransform}
\mathcal{G}_{a_{i},b_{i}}=\{g_{a_{i},b_{i}}(x)=a_{i}x_{ij}+b_{i},~j=1,\ldots,n_{i},~i=1,2\}.
\end{eqnarray}
Under the group $\mathcal{G}$ the estimation problem is invariant and  the form of the affine equivariant estimators is given as 
\begin{eqnarray}\label{1c}
\delta_{1c}(\underline{X_{(1)}},\underline{T})=X_{1(1)}-cT_1,
\end{eqnarray}
where $c$ is a constant. The BAEE of $\mu_{1}$ is given by the following theorem.
\begin{theorem}
Let  $U_1 = \frac{X_{1(1)} - \mu_1}{\sigma_1}$, which
follow an exponential $Exp\left(\frac{1}{n_1}\right)$ distribution
and $V_1=\frac{T_1}{\sigma_1}$, which follow a
$Gamma(n_1-1,1)$ distribution. Then, under an arbitrary bowl-shaped loss function $L(t)$, the
(unrestricted) BAEE
$\delta_{01}(\underline{X_{(1)}},\underline{T})=X_{1(1)}-c_{01}T_{1}$,
where $c_{01}$ is the unique solution of the equation
\begin{eqnarray}\label{c01}
E\left(L^{\prime}\left(U_1-c_{01}V_1\right)V_1\right)=0
\end{eqnarray}
\end{theorem}
\noindent{\bf Proof.} Proof of this theorem is straight forward 
\hfill\(\Box\)
\begin{example}\rm
Consider
$L_1=\left(\frac{\delta-\mu_1}{\sigma_1}\right)^2$, the squared error loss. Under this loss
function $c_{01}$ is the unique solution of the equation
\begin{eqnarray*}
\int_{0}^{\infty}\int_{0}^{\infty}\left(u_1-c_{01}v_1\right)v_1e^{-n_1u_1}e^{-v_1}v_1^{n_1-2} du_1 dv_1=0.
\end{eqnarray*}
So we have   $c_{01}=\frac{1}{n_1^2}$ and the BAEE of $\mu_1$ with respect to $L_1$ is obtained as
\begin{eqnarray*}
\delta_{01}^Q=X_{1(1)}-\frac{T_{1}}{n^2_1}.
\end{eqnarray*}
\end{example}
\begin{example}\rm
    Consider $L_2=e^{a\left(\frac{\delta-\mu_1}{\sigma_1}\right)}-a\left(\frac{\delta-\mu_1}{\sigma_1}\right)-1$, the linex loss function. Under this loss function $c_{01}$ is the unique solution
    \begin{eqnarray*}
       \int_{0}^{\infty}\int_{0}^{\infty} a (e^{a(u_1-c_{01}v_1)} - 1) v_1e^{-n_1u_1}e^{-v_1}v_1^{n_1-2}   du_1 dv_1=0.
    \end{eqnarray*}
We get   $c_{01}=\frac{1}{a}\left(\left(\frac{n_1}{n_1-a}\right)^{1/n_1} -1\right)$, $n_1>a$. So BAEE of $\mu_1$ is
 $$\delta_{01}^L=X_{1(1)} - \frac{1}{a}\left(\left(\frac{n_1}{n_1-a}\right)^{1/n_1} - 1\right) T_{1}.$$
\end{example}
\begin{remark}\rm
 We note that the UMVUE of $\mu_1$ is
    $$\delta_{1MV}=X_{1(1)}-\frac{T_1}{n_1(n_1-1)}$$
and the maximum likelihood estimator (MLE) of $\mu_1$ is
$\delta_{1ML}=X_{1(1)}$. Both these estimators belong to the class
of affine estimators in (\ref{1c}), so under quadratic and linex
losses, these two estimators are inadmissible.
\end{remark}

\subsection{Improvement over the BAEE of $\mu_1$}\label{sec2}
Now we will show that BAEE of $\mu_1$ is inadmissibly. To find the dominating estimator take an estimator of the form (see \cite{vijayasree1995componentwise})
\begin{eqnarray}\label{equiphi1}
\delta_{\phi_{1}}\left(\underline{X_{(1)}},\underline{T}\right) =X_{1(1)}-T_{1}\phi_{1}\left(\frac{T_{2}}{T_{1}}\right)=X_{1(1)}-T_{1}\phi_{1}\left(W\right),
\end{eqnarray}
where $W=\frac{T_2}{T_1}$ and $\phi_{1}$ is a positive measurable function. Note that the BAEE belongs to this class.

In the following theorem we prove that the BAEE of $\mu_1$,
$\delta_{01}$, is inadmissible.
\begin{theorem}\label{thmst1}
Let $b_{01}$ be the unique solution of
\begin{eqnarray}\label{b01}
E\left(L^{\prime}(U_1-b_{01}Z)\right)=0,
\end{eqnarray}
where $Z$ follow a $Gamma(n_1+n_2-1,1)$ distribution and $U_1$
follow an exponential $Exp(1/n_1)$ distribution.
Consider an estimator of the form
\begin{eqnarray}
\delta_{1S}=X_{1(1)}-\phi_{01}(W)T_1~~~ \mbox{with}~~~ \phi_{01}(w)=\min\left\{c_{01},b_{01}(1+w)\right\},
\end{eqnarray}
where $c_{01}$ is given in (\ref{c01}). Then the risk of the
estimator $\delta_{1S}$ is nowhere larger than that of $\delta_{01}$
under a bowl shaped location-scale invariant loss function $L(t)$  provided
\begin{eqnarray}\label{sttheq1}
E\left(L^{\prime}\left(U_1-b_{01}V_1\right)V_1\right) \ge 0.
\end{eqnarray}
\end{theorem}
\noindent \textbf{Proof:} The risk function of the estimator $\delta_{\phi_{1}}\left(\underline{X_{(1)}},\underline{T}\right)$ is
\begin{eqnarray*}
R(\delta_{\phi_{1}})=E\left(L\left(\frac{X_{1(1)}-T_{1}\phi_{1}(W)-\mu_1}{\sigma_1}\right)\right)
=E^WE\left(L\left(\frac{X_{1(1)}-T_{1}\phi_{1}(W)-\mu_1}{\sigma_1}\right)\bigg|W=w\right).
\end{eqnarray*}
Define
\begin{eqnarray*}
R_1(c)=E\left(L\left(\frac{X_{1(1)}-cT_{1}-\mu_1}{\sigma_1}\right)
\bigg|W=w\right).
\end{eqnarray*}
It can be seen that the conditional distribution of
$V_1=T_1/\sigma_1$ given $W=w$ follow a gamma
$Gamma(n_1+n_2-2,(1+\eta w)^{-1})$ distribution, where
$\eta=\frac{\sigma_1}{\sigma_2}\le1$. Let, $U_1 =
\dfrac{X_{1(1)} - \mu_1}{\sigma_1}$. Hence we have
\begin{eqnarray*}
R_1(c)=E_{\eta}\left(L\left(U_1-V_1c\right)\big|W=w\right).
\end{eqnarray*}
The conditional risk function $R_1(\sigma_1,\sigma_2,c)$ is
minimized at $c_{\eta}(w)$, where $c_{\eta}(w)$ is the unique
solution of
\begin{eqnarray}\label{st_2}
E_{\eta}\left(L^{\prime}\left(U_1-V_1c_{\eta}(w)\right)\big|W=w\right)=0,
\end{eqnarray}
where in the equation (\ref{st_2}) the expectation has been taken with respect to the density $\zeta_{\eta}(u_1,v_1|w) \propto v_1^{n_1+n_2-2} e^{- n_1 u_1 -
v_1 (1+ \eta w)}$.
It can be easily seen that
$\frac{\zeta_{\eta}(u_1,v_1|w)}{\zeta_{1}(u_1,v_1|w)}$ is non
decreasing in $v_1$.  Also we have $L^{\prime}(u_1-v_1c_{\eta}(w))$ is decreasing in $v_1$ as $c_{\eta}(w)>0$. $c_{\eta}(w)>0$ is guaranteed  as  $\phi(w)$ is a positive function. 
Consequently, using Lemma 3.4.2. in
\cite{lehmann2005testing}, we have
$$E_{\eta}\left(L^{\prime}\left(U_1-V_1c_{1}(w)\right)\big|W=w\right)
\le
E_{1}\left(L^{\prime}\left(U_1-V_1c_{1}(w)\right)\big|W=w\right)
= 0.$$
So from (\ref{st_2}) we have 
\begin{equation} \label{st_3}
	E_{\eta}\left(L^{\prime}\left(U_1-V_1c_{1}(w)\right)\big|W=w\right)
\le
E_{\eta}\left(L^{\prime}\left(U_1-V_1c_{\eta}(w)\right)\big|W=w\right)
= 0.
\end{equation}
 As we have $E_{\eta}L^{\prime}(U_1-V_1c_{\eta}(w))|W=w)$ is decreasing in $c_{\eta}(w)$. The inequality (\ref{st_3}) implies that
 \begin{eqnarray}
c_{\eta}(w) \le c_{1}(w).
\end{eqnarray}

Again $c_{1}(w)$ is the unique solution of
\begin{eqnarray}
\int_{0}^{\infty}\int_{0}^{\infty}L^{\prime}\left(u_1-v_1c_{1}(w)\right)e^{-n_1
u_1}v^{n_1+n_2-2}_1 e^{-v_1(1+w)}dv_1du_1=0.
\end{eqnarray}
By the transformation $u_1=u_1$, $z=v_1(1+w)$ we get
\begin{eqnarray}\label{sttheq4}
\int_{0}^{\infty}\int_{0}^{\infty}L^{\prime}\left(u_1-\frac{z
c_1(w)}{1+w}\right) e^{- n_1
u_1}z^{n_1+n_2-2}e^{-z}dzdu_1=0.
\end{eqnarray}
Comparing (\ref{sttheq4}) with (\ref{b01}), we get $c_1(w)=b_{01}(1+w)$. Moreover (\ref{c01}) and (\ref{sttheq1}) ensure that
\begin{eqnarray}
 b_{01}\le c_{01} .
\end{eqnarray}
Consider a function $\phi_{01}(w)=\min\left\{c_{01},b_{01}(1+w)\right\}$. Then we have
$c_{\eta}(w) <\phi_{01}(w)<c_{01}$ on a set of positive probability. Since
$R_1(c)$ is strictly bowl shaped, so for $c>c_{\eta}(w)$,
$R_1(c)$ is an increasing function of $c$. Therefore,
$R(\delta_{1S}) \le R(\delta_{01})$. This completes the proof of the theorem. \hfill $\blacksquare$\\
\begin{corollary}\label{coro1}
Consider an estimator of the form (\ref{equiphi1}), that is
$\delta_{\phi_{1}}\left(\underline{X_{(1)}},\underline{T}\right)
=X_{1(1)}-T_{1}\phi_{1}\left(W\right).$ Define
$\phi_1^{*}(w)=\min\{\phi_{1}\left(W\right),b_{01}(1+w)\}$. Then the
estimator
$$\delta_{\phi_1^{*}}=X_{1(1)}-T_{1}\phi^{*}_{1}\left(W\right)$$
improve upon $\delta_{\phi_{1}}$ provided
$P(\phi_{1}(W)>b_{01}(1+W)) \ne 0$ under a bowl-shaped location
invariant loss function $L(t)$.
\end{corollary}

%

\begin{example}\rm
For  $L_1=\left(\frac{\delta-\mu_1}{\sigma_1}\right)^2$. Then the  estimator $$\delta^Q_{1S}=X_{1(1)}-\min\left\{\frac{1}{n_1^2},\frac{1}{n_1(n_1+n_2-1)}(1+W)\right\} T_1$$ improve upon the BAEE of $\mu_1$.  It is pointed out that the improved estimator for $\mu_1$, $\delta^Q_{1S}$, is the same as the estimator discussed in Corollary 3.3 by \cite{vijayasree1995componentwise} for $k=2$. The estimator which improve upon the UMVUE is obtained as $$\delta^{Q}_{1MV}=X_{1(1)}-\min\left\{\frac{1}{n_1(n_1-1)},\frac{1}{n_1(n_1+n_2-1)}(1+W)\right\}T_1$$

\end{example}
\begin{example}
For $L_2=e^{a\left(\frac{\delta-\mu_1}{\sigma_1}\right)}-a\left(\frac{\delta-\mu_1}{\sigma_1}\right)-1$ and we have for $n_1>a$,   $$b_{01}=\frac{1}{a}\left[\left(\frac{n_1}{n_1-a}\right)^{1/(n_1+n_2-1)} - 1\right]~~~\mbox{ and }~~~c_{01}=\frac{1}{a}\left[\left(\frac{n_1}{n_1-a}\right)^{1/n_1} - 1\right]$$ 
Then the  estimator
$$\delta^L_{1S}=X_{1(1)}-\min\left\{c_{01},b_{01}(1+W)\right\} T_1$$
improve upon the BAEE of $\mu_1$. 

\noindent The estimator which improve upon the UMVUE of $\mu_1$ is obtained as
$$\delta^{L}_{1MV}=X_{1(1)}-\min\left\{\frac{1}{n_1(n_1-1)},b_{01}(1+W)\right\} T_1$$
\end{example}


Next, we consider another class of estimators of the form

\begin{eqnarray}\label{st2}
\delta_{\phi_{1}}^*\left(\underline{X_{(1)}},\underline{T}\right)=X_{1(1)}-T_{1}\phi_{1}^*\left(W,W_1\right),
\end{eqnarray}
where $W=\frac{T_2}{T_1}$, $W_1=\frac{X_{2(1)}}{T_1}$ and
$\phi_{1}^*$ is a positive measurable function.
\begin{theorem}\label{thmst2}
Let $b_{01}^*$ be the unique solution of
\begin{eqnarray}\label{b02}
E\left(L^{\prime}(U_1-b_{01}^*Z_1)\right)=0,
\end{eqnarray}
where $Z_1$ follow a $Gamma(n_1+n_2,1)$ distribution and $U_1$
follow a standard exponential $Exp(1/n_1)$ distribution. Consider an
estimator of the form
\begin{eqnarray}
\delta_{1S}^*=X_{1(1)}-\phi_{01}^*(W, W_1)T_1
\end{eqnarray}
with
\begin{eqnarray}
\phi_{01}^*(w,w_1)=\left\{\begin{array}{lcl}
\min\left\{c_{01},b_{01}^*(1+w+w_1)\right\}, ~~ w_1 > 0 \\
c_{01},~~~~~~~~~~~~~~~~~~~~~~~~~~~~\text{otherwise}
\end{array}
\right.,
\end{eqnarray}
where $c_{01}$ is given in (\ref{c01}). Then the risk of the
estimator $\delta_{1S}^*$ is nowhere larger than that of
$\delta_{01}$ under a bowl shaped scale invariant loss function $L(t)$
provided
\begin{eqnarray}\label{sttheq11}
E\left(L^{\prime}\left(U_1-b_{01}^*V_1\right)V_1\right) \ge 0.
\end{eqnarray}
\end{theorem}
\noindent \textbf{Proof:} The proof is analogous to Theorem
\ref{thmst1} and that's why is omitted.
\hfill $\blacksquare$\\
\begin{corollary}\label{coro6}
Consider an estimator of the form (\ref{st2}), that is $\delta_{\phi_{1}}^*\left(\underline{X_{(1)}},\underline{T}\right)=X_{1(1)}-T_{1}\phi_{1}^*\left(W,W_1\right),$

Define $$\psi_1^{*}(w,w_1)=\left\{\begin{array}{lcl}
        \min\{\phi_{1}^*\left(w,w_1\right),b^*_{01}(1+w+w_1)\}, ~~~w_1>0\\\\

        \phi_{1}^*\left(w,w_1\right),~~~~~~~~~~~~~~otherwise
    \end{array}
\right.,$$
Then the estimator
        $$\delta_{\psi_1^{*}}=X_{1(1)}-T_{1}\psi^{*}_{1}\left(W,W_1\right)$$
        improve upon $\delta_{\phi_{1}}$ provided $P(\phi_{1}^*\left(w,w_1\right)>b^*_{01}(1+w+w_1)) \ne 0$ under a bowl-shaped location invariant loss function $L(t)$.
    \end{corollary}

 It can be easily seen that the UMVUE and MLE of $\mu_1$ belong to the class of estimators in (\ref{st2}). As a consequence of Corollary \ref{coro6} we get the following dominance result.

\begin{corollary}\label{corm1}
With respect to a bowl shaped location invariant loss function $L_1\left(\frac{\delta-\mu_1}{\sigma_1}\right)$,
         the estimator $\delta_{1MV}$ is inadmissible and dominated by
        $$\delta^{I2}_{1MV}=\left\{\begin{array}{lcl}
        X_{1(1)}-\min\left\{\frac{1}{n_1(n_1-1)},b^*_{01}(1+W+W_1)\right\} T_1, W_1>0\\\\
        X_{1(1)}-\frac{1}{n_1(n_1-1)}  T_1, ~~~~~~~~~~~~~~~~~~~~~~~~~otherwise
            \end{array}
    \right.,
        $$
\end{corollary}

\begin{example}
    Under the squared error loss $L_1=\left(\frac{\delta-\mu_1}{\sigma_1}\right)^2$function the
    estimator of the form
    \begin{eqnarray}
        \delta_{1S}^{*,Q}=X_{1(1)}-\phi_{01}^*(W, W_1)T_1
    \end{eqnarray}
    with
    \begin{eqnarray}
        \phi_{01}^*(w,w_1)=\left\{\begin{array}{lcl}
            \min\left\{\frac{1}{n_1^2},\frac{(1+w+w_1)}{n_1(n_1+n_2)}\right\}, ~~ w_1 > 0 \\
            \frac{1}{n_1^2},~~~~~~~~~~~~~~~~~~~~~~~~~~~~\text{otherwise}
        \end{array}
        \right.,
    \end{eqnarray}
improves upon the BAEE of $\mu_1$. The UMVUE of $\mu_1$ is inadmissible and dominated by
$$\delta^{2,Q}_{1MV}=\left\{\begin{array}{lcl}
    X_{1(1)}-\min\left\{\frac{1}{n_1(n_1-1)},\frac{(1+w+w_1)}{n_1(n_1+n_2)}\right\} T_1, W_1>0\\\\
    X_{1(1)}-\frac{1}{n_1(n_1-1)} T_1, ~~~~~~~~~~~~~~~~~~~~~~~~~otherwise
\end{array}
\right.,
$$
\end{example}
\begin{example}
For $L_2=e^{a\left(\frac{\delta-\mu_1}{\sigma_1}\right)}-a\left(\frac{\delta-\mu_1}{\sigma_1}\right)-1$. Under this loss for $n_1>a$ function we have
$$b^*_{01}=\frac{1}{a}\left[\left(\frac{n_1}{n_1-a} \right)^{1/(n_1+n_2)} - 1\right]~~~~\mbox{ and }~~~c_{01}=\frac{1}{a}\left[\left(\frac{n_1}{n_1-a}\right)^{1/n_1} - 1 \right]$$
    The improved estimators for $\mu_1$ can be easily derived form Corollary \ref{corm1} by substituting the values of $b^*_{01}$ and $c_{01}$.
\end{example}

\subsection{A class of improved estimator for $\mu_1$}
In this section, we derive a class of improved estimators. For this
we consider a class of estimators of the form
\begin{eqnarray}\label{1kubo}
\delta_{\varphi_1}=X_{1(1)}-T_1\varphi_1(W).
\end{eqnarray}

 Let $U_1 = \dfrac{X_{1(1)} - \mu_1}{\sigma_1}$ with
pdf $h(u_1) = n_1 e^{-n_1 u_1}$, $u_1 > 0$ and $V_1 =
\dfrac{T_1}{\sigma_1}$.  The joint pdf of $V_1$, $W$ is provided as
\begin{eqnarray*}
g_{\eta}(v_1,w) \propto \eta^{n_2-1}e^{-v_1\left(1+\eta
w\right)}w^{n_2-2}v_1^{n_1+n_2-3},~~v_1>0,w>0,0<\eta \le 1
\end{eqnarray*}
and 
 we define the joint cdf of $V_1$, $W$ as, $$G_{\eta}(v_1, z) = \int_0^{z} g_{\eta}(v_1,x) dx$$.

\begin{theorem}\label{thmkubo1}
Let the function $\varphi_{1}(z)$ satisfies the following
assumptions
\begin{itemize}
\item[(i)] $\varphi_{1}(z)$ is non-decreasing and
$\lim_{z\rightarrow\infty}\varphi_{1}(z)=c_{01}$,
\item[(ii)]$\int_{0}^{\infty}
\int_{0}^{\infty}L^{\prime}\left(u_1-v_1\varphi_{1}(z)\right)h(u_1)v_1G_{1}(v_1,z)du_1dv_1 \le0.$
\end{itemize}
Then the risk of the estimator $\delta_{\varphi_{1}}$ given in (\ref{1kubo}) dominates $\delta_{01}$ under a bowl shaped location-scale invariant loss function $L(t)$.
\end{theorem}
\noindent{\textbf{Proof}:} The risk difference of $\delta_{01}$ and $\delta^{1}_{\varphi_{1}}$ can be written as
\begin{eqnarray}\nonumber
R_{\delta_{01}}-R_{\delta^{1}_{\varphi_{1}}}&=&
EL\left(\frac{X_{1(1)}-T_1c_{01}-\mu_1}{\sigma_1}\right)-EL\left(\frac{X_{1(1)}-T_1\varphi_{1}(W)-\mu_1}{\sigma_1}\right)\\ \nonumber
&=&-E\int_{1}^{\infty}L^{\prime}\left(U_1-V_1\varphi_{1}(tW)\right)\varphi_{1}^{\prime}(tW)V_1Wdt\\ \nonumber
&=&-\int_{0}^{\infty}\int_{0}^{\infty}\int_{0}^{\infty}\int_{1}^{\infty}L^{\prime}\left(u_1-v_1\varphi_{1}(tw)\right)
\varphi_{1}^{\prime}(tw)v_1wh(u_1)g_{\eta}(v_1,w)dtdwdv_1du_1,
\end{eqnarray}
Now we use the transformation $wt=z$ and $x=z/t$, then the risk
difference reduces to
\begin{eqnarray}\nonumber
R_{\delta_{01}}-R_{\delta^{1}_{\varphi_{1}}}&=&-\int_{0}^{\infty}\int_{0}^{\infty}\int_{0}^{\infty}
\int_{0}^{z}L^{\prime}\left(u_1-v_1\varphi_{1}(z)\right)
\varphi_{1}^{\prime}(z)v_1h(u_1)g_{\eta}(v_1,x)dxdzdv_1du_1\\\nonumber
&=&-\int_{0}^{\infty}\int_{0}^{\infty}\int_{0}^{\infty}L^{\prime}\left(u_1-v_1\varphi_{1}(z)\right)
\varphi_{1}^{\prime}(z)v_1h(u_1)G_{\eta}(v_1,z)dzdv_1du_1\\\nonumber
&=&-\int_{0}^{\infty}\int_{0}^{\infty}\varphi_{1}^{\prime}(z)\left[\int_{0}^{\infty}
L^{\prime}\left(u_1-v_1\varphi_{1}(z)\right)h(u_1)du_1\right]
v_1G_{\eta}(v_1,z)dv_1dz\\\nonumber
&=&-\int_{0}^{\infty}\int_{0}^{\infty}\varphi_{1}^{\prime}(z)\left[\int_{0}^{\infty}
L^{\prime}\left(u_1-v_1\varphi_{1}(z)\right)h(u_1)du_1\right]
v_1\frac{G_{\eta}(v_1,z)}{G_{1}(v_1,z)}G_{1}(v_1,z)dv_1dz
\end{eqnarray}

Take $\mathcal{F}(v_1)=-\int_{0}^{\infty}
L^{\prime}\left(u_1-v_1\varphi_{1}(z)\right)h(u_1)du_1$. $L(t)$
is strictly bowl shaped and non monotone function in $\mathbb{R}$.
So we have $\mathcal{F}(v_1)<0$ for $v_1<v_1^*$ and
$\mathcal{F}(v_1)>0$ for $v_1>v_1^*$ where $v_1^*$ is the solution
of $\mathcal{F}(v_1^*)=0$. Also it can be easily seen that
$\frac{g_{\eta}\left(x,v_1\right)}{g_{1}\left(x,v_1\right)}$ is
non-decreasing in $x$ for $x>0$ and consequently,
$\frac{G_{\eta}\left(x,v_1\right)}{G_{1}\left(x,v_1\right)}$ is
non-decreasing in $x$ for $x>0$.Then by Lemma (2.1) of
\cite{kubokawa1994unified} we have
\begin{eqnarray}\nonumber
R_{\delta_{01}}-R_{\delta^{1}_{\varphi_{1}}}
&=&-\int_{0}^{\infty}\int_{0}^{\infty}\varphi_{1}^{\prime}(z)\left[\int_{0}^{\infty}
L^{\prime}\left(u_1-v_1\varphi_{1}(z)\right)h(u_1)du_1\right]
v_1\frac{G_{\eta}(v_1,z)}{G_{1}(v_1,z)}G_{1}(v_1,z)dv_1dz\\\nonumber
&\ge &-\int_{0}^{\infty}\varphi_{1}^{\prime}(z)\frac{G_{\eta}(v_1^*,z)}{G_{1}(v_1^*,z)}\int_{0}^{\infty}
\int_{0}^{\infty}L^{\prime}\left(u_1 - v_1\varphi_{1}(z)\right)h(u_1)du_1
v_1G_{1}(v_1,z)dv_1dz\\\nonumber
&=&-\int_{0}^{\infty}\varphi_{1}^{\prime}(z)\frac{G_{\eta}(v_1^*,z)}{G_{1}(v_1^*,z)}\left[\int_{0}^{\infty}
\int_{0}^{\infty}L^{\prime}\left(u_1 - v_1\varphi_{1}(z)\right)h(u_1)v_1G_{1}(v_1,z)du_1
dv_1\right]dz.
\end{eqnarray}
Hence by the condition $(ii)$ the risk difference in non negative. This proves the result.
	\begin{remark} \label{rembz1}
		The limit case in Theorem \ref{thmkubo1} provides the \cite{brewster1974improving} - type estimator (see \cite{kubokawa1994unified} for details). 
		That means, $\delta_{\varphi_{BZ}}=X_{1(1)}-T_1\varphi_{BZ}(W)$ is the \cite{brewster1974improving}- type estimator of $\mu_1$, where $\varphi_{BZ}$ is given in Equation \eqref{eqbz}
		\begin{equation} \label{eqbz}
			\int_{0}^{\infty} \int_{0}^{\infty}L^{\prime}\left(u_1-v_1\varphi_{BZ}(z)\right)h(u_1)v_1G_{1}(v_1,z)du_1dv_1 = 0.
		\end{equation}
		This type of estimator is, usually, a generalized Bayes estimator. Indeed, in our case, we can prove that $\delta_{\varphi_{BZ}}$ is a generalized Bayes estimator of $\mu_1$ under the improper prior $\pi(\mu_1, \sigma_1, \sigma_2) = \frac{1}{\sigma_1 \sigma_2}$, $\mu_1 \in R$, $0< \sigma_1 < \sigma_2$ and the bowl-shaped loss function $L(t)$.
		
		It can be shown (see \cite{marchard2005}) that this Bayes estimator is of the form
		$$ \delta = \delta_0 + T_1 g_0(W) = X_{1(1)}+T_1 g(W) .$$ An estimator of the form $X_{1(1)}+T_1 g(W)$ is the generalized Bayes with respect to the prior $\pi$ provided the corresponding posterior expected loss
		
		$$E \Big( L\left(\frac{X_{1(1)}+T_1 g(W) - \mu_1}{\sigma_1} \right) \Big| X_{1(1)}=x_1, T_1=t_1, W=w \Big)$$ is minimum for this estimator.
		
		Equivalently
		
		$$\int_{0}^{\infty} \int_{0}^{\sigma_2} \int_{-\infty}^{x_1} L\Big( \frac{x_1+t_1 g(w) - \mu_1}{\sigma_1} \Big) \frac{1}{\sigma_1 \sigma_2} \frac{n_1}{\sigma_1} e^{-n_1 \frac{x_1-\mu_1}{\sigma_1}}\,
		\frac {1}{\sigma_1^{n_1-1} \sigma_2^{n_2-1}} \frac{1}{ \Gamma(n_1-1) \Gamma(n_2-1)}\, t_1^{n_1+n_2-3} \times$$ $$ \times e^{-t_1/\sigma_1}\, e^{-t_1 w/\sigma_2} w^{n_2-2} \, d\mu_1 d\sigma_1 d\sigma_2 $$
		should be minimized in $g$. Using the transformation $U_1 = \frac{X_{1(1)} - \mu_1}{\sigma_1}, V_1 = \frac{T_1}{\sigma_1}$ and $\eta = \frac{\sigma_1}{\sigma_2}$ 
		we see that Jacobean $J\left(\frac{(\mu_1, \sigma_1, \sigma_2)}{(u_1, v_1, \eta)}\right) = \frac{t_1^3}{v_1^4 \eta^2}$ and the above expression can be rewritten as
		$$\int_{0}^{\infty} \int_{0}^{\infty} \int_{0}^{1} L(u_1+v_1 g(w) ) \frac{1}{ \Gamma(n_1-1) \Gamma(n_2-1)} e^{-n_1 u_1}\, v_1^{n_1+n2-3} \frac{1}{t_1} \eta^{n_2-2} w^{n_2-2} e^{- v_1 (1+\eta w)} \, d\eta du_1 dv_1 $$.
		
		Minimizing in $g$ the latter expression and letting $x = \eta w$, we conclude that $g$ is such that
		\begin{equation} \label{eqg}
			\int_{0}^{\infty} \int_{0}^{\infty} L'(u_1+v_1 g(w) ) v_1 n_1 e^{-n u_1} v_1^{n_1+n_2-3} \int_{0}^{w} x^{n_2-2} e^{-v_1(1+x)} dx du_1 dv_1 = 0
		\end{equation}
		So for the choice of $g(w) = -\varphi_{BZ}(w)$, equation \eqref{eqg} is the same as \eqref{eqbz} and the estimator $\delta_{\varphi_{BZ}}$ is the generalized Bayes estimator of $\mu_1$ under the $L(t)$ loss function and the prior $\pi$.
	\end{remark}
\begin{example}
Suppose the function $\varphi_{1}(z)$ satisfies the following assumptions
\begin{itemize}
    \item[(i)] $\varphi_{1}(z)$ is non-decreasing and
   $\lim_{z\rightarrow\infty}\varphi_{1}(z)=  \frac{1}{n_1^2}$,
  \item[(ii)]$\varphi_{1}(z)  \ge \varphi_{1BZ}(z)$,
    with
    $$\varphi_{1BZ}(z)=\frac{\int_{0}^{\infty}\int_{0}^{z}e^{-v_1(1+x)}x^{n_2-2}v_1^{n_1+n_2-2}dxdv_1}{n_1\int_{0}^{\infty}\int_{0}^{z}e^{-v_1(1+x)}x^{n_2-2}v_1^{n_1+n_2-1}dxdv_1}$$
\end{itemize}
Then the risk of the estimator $\delta_{\varphi_{1}}$ given in (\ref{1kubo}) dominates $\delta_{01}$ under squared error loss $L_1\left(\frac{\delta-\mu_1}{\sigma_1}\right)$.
\end{example}

 \begin{example}
Suppose the function $\varphi_{1}(z)$ satisfies the following assumptions
\begin{itemize}
    \item[(i)] $\varphi_{1}(z)$ is non-decreasing and
   $\lim_{z\rightarrow\infty}\varphi_{1}(z)= \frac{1}{a}\left[\left(\frac{n_1}{n_1-a}\right)^{1/n_1} - 1\right]$,
   \item[(ii)] $\int_{0}^{\infty} \int_{0}^{\infty} a \left(e^{a (u_1-v_1\varphi_{1}(z) )} -1\right)h(u_1)v_1G_{1}(v_1,z)du_1 dv_1 \le 0.$
\end{itemize}
Then the risk of the estimator $\delta_{\varphi_{1}}$ given in (\ref{1kubo}) dominates $\delta_{01}$ under linex loss $L_2\left(\frac{\delta-\mu_1}{\sigma_1}\right)$.
\end{example}

\section{Improved estimation of $\mu_2$}\label{sec3}
Now, we consider the problem of estimating the location parameter
$\mu_{2}$ knowing that the order restriction
$\sigma_{1}\leq\sigma_{2}$ exists. Considering the affine group of
transformations (\ref{eqtransform}), see Section 2, the form of the
affine equivariant estimators can be obtained as
\begin{eqnarray}\label{2c}
\delta_{2c}(\underline{X_{(1)}},\underline{T})=X_{2(1)} - cT_2,
\end{eqnarray}
where $c$ is a constant. The following theorem provides the BAEE of
$\mu_{2}$.
\begin{theorem}
Let $U_2=\frac{X_{2(1)} - \mu_2}{\sigma_2}$, which follows an
exponential $Exp\left(\frac{1}{n_2}\right)$ distribution and $V_2 =
\frac{T_2}{\sigma_2}$, which follows a $Gamma(n_2-1,1)$
distribution. Then under a arbitrary bowl-shaped loss function $L(t)$ the
(unrestricted) best affine equivariant estimator of $\mu_{2}$ is
$\delta_{02}(\underline{X_{(1)}},\underline{T})=X_{2(1)} - c_{02}T_{2}$,
where $c_{02}$ is the unique solution of the equation
\begin{eqnarray}\label{c02}
E\left(L^{\prime}\left(U_2 - c_{02}V_2\right)V_2\right)=0
\end{eqnarray}
\end{theorem}
\noindent{\bf Proof.} The proof is simple and hence omitted for the
sake of brevity. \hfill\(\Box\)

\begin{example}\rm Consider the squared error loss
function $L_1=\left(\frac{\delta-\mu_2}{\sigma_2}\right)^2$. Under
this loss function $c_{02}$ is the unique solution of the equation
\begin{eqnarray} \label{eqc2quad}
\int_{0}^{\infty}\int_{0}^{\infty}\left(u_2 - c_{02}v_2\right)v_2e^{-n_2u_2}e^{-v_2}v_2^{n_2-2}
du_2 dv_2=0.
\end{eqnarray}
After some calculation, from (\ref{eqc2quad}), we get that $c_{02}=\frac{1}{n_2^2}$ and the BAEE of $\mu_2$ with respect to
$L_1$ is obtained as
\begin{eqnarray*}
\delta_{02}^Q=X_{2(1)}-\frac{T_{2}}{n_2^2}.
\end{eqnarray*}
\end{example}
\begin{example}\rm
Consider the linex loss function
$L_2=e^{a\left(\frac{\delta-\mu_2}{\sigma_2}\right)}-a\left(\frac{\delta-\mu_2}{\sigma_2}\right)-1$.
Under this loss function $c_{02}$ is the unique solution
\begin{eqnarray} \label{eqc2linex}
\int_{0}^{\infty}\int_{0}^{\infty} a (e^{a(u_2 - c_{02}v_2)} - 1)
v_2e^{-n_2u_2}e^{-v_2}v_2^{n_2-2} du_2 dv_2=0.
\end{eqnarray}
From equation (\ref{eqc2linex}), we get that $c_{02}=\frac{1}{a}\left(\left(\frac{n_2}{n_2-a} \right)^{1/n_2} -1\right)$,
$n_2>a$. So, the BAEE of $\mu_2$ is $$\delta_{02}^L=X_{2(1)}  - \frac{1}{a}\left(\left(\frac{n_2}{n_2-a}\right)^{1/n_2} - 1\right)
T_{2}.$$
\end{example}
\begin{remark}\rm
We note that the UMVUE of $\mu_2$ is
$$\delta_{2MV}=X_{2(1)}-\frac{T_2}{n_2(n_2-1)}$$ and the maximum
likelihood estimator (MLE) of $\mu_2$ is $\delta_{2ML}=X_{2(1)}$.
Both these estimators belong to the class of affine estimators in
(\ref{2c}), so under quadratic and linex losses, these two
estimators are inadmissible.
\end{remark}

\subsection{Improvement over the BAEE of $\mu_2$}
In this subsection, we want to prove that the BAEE of $\mu_2$ is
inadmissible by deriving an improved estimator of Stein-type.
Consider an estimator of the form
\begin{eqnarray}\label{equiphi2}
\delta_{\phi_{2}}\left(\underline{X_{(1)}},\underline{T}\right)
=X_{2(1)} - T_{2}\phi_{2}\left(\frac{T_{1}}{T_{2}}\right)=X_{2(1)} - T_{2}\phi_{2}\left(W^*\right),
\end{eqnarray}
where $W^*=\frac{T_1}{T_2}$ and $\phi_{2}$ is a positive measurable function.
In the following theorem we prove that the BAEE of $\mu_2$,
$\delta_{02}$, is inadmissible.
\begin{theorem}\label{thmst31}
Let $b_{02}$ be the unique solution of
\begin{eqnarray}\label{bstar01}
E\left(L^{\prime}(U_2 - b_{02}Z^*)\right)=0,
\end{eqnarray}
where $Z^*$ follow a $Gamma(n_1+n_2-1,1)$ distribution and $U_2$
follow a standard exponential $Exp(1/n_2)$ distribution. Consider an
estimator of the form
\begin{eqnarray}
\delta_{2S}=X_{2(1)} - \phi_{02}(W^*)T_2~~~ \mbox{with}~~~
\phi_{02}(w^*)=\max\left\{c_{02},b_{02}(1+w^*)\right\},
\end{eqnarray}
where $c_{02}$ is given in (\ref{c02}). Then the risk of the
estimator $\delta_{2S}$ is nowhere larger than that of $\delta_{02}$
under a bowl shaped invariant loss function $L(t)$ provided
\begin{eqnarray}\label{sttheqbstar01}
E\left(L^{\prime}\left(U_2 - b_{02}V_2\right)V_2\right)  \le 0.
\end{eqnarray}
\end{theorem}
\noindent \textbf{Proof:} The risk function of the estimator
$\delta_{\phi_{2}}\left(\underline{X_{(1)}},\underline{T}\right)$ is
\begin{eqnarray*}
R(\delta_{\phi_{2}})=E\left(L\left(\frac{X_{2(1)} - T_{2}\phi_{2}(W^*)-\mu_2}{\sigma_2}\right)\right)
=E^WE\left(L\left(\frac{X_{2(1)} - T_{2}\phi_{2}(W^*)-\mu_2}{\sigma_2}\right)\bigg|W^*=w^*\right).
\end{eqnarray*}
Define
\begin{eqnarray*}
R_2(c)=E\left(L\left(\frac{X_{2(1)} - cT_{2}-\mu_2}{\sigma_2}\right)
\bigg|W^*=w^*\right).
\end{eqnarray*}
It can be seen that the conditional distribution of
$V_2=T_2/\sigma_2$ given $W^*=w^*$ follow a
$Gamma(n_1+n_2-2,(1+w^*/\eta)^{-1})$ distribution, where
$\eta=\frac{\sigma_1}{\sigma_2}\le1$. Let $U_2 = \dfrac{X_{2(1)} -
\mu_2}{\sigma_2}$ Hence we have
\begin{eqnarray*}
R_2(c)=E_{\eta}\left(L\left(U_2 - V_2c\right)\big|W^*=w^*\right).
\end{eqnarray*}
The conditional risk function $R_2(\sigma_1,\sigma_2,c)$ is
minimized at $c^*_{\eta}(w^*)$, where $c^*_{\eta}(w^*)$ is the
unique solution of

\begin{eqnarray} \label{eq3_9}
E_{\eta}\left(L^{\prime}\left(U_2 - V_2c^*_{\eta}(w^*)\right) \big|W^*=w^*\right)=0.
\end{eqnarray}
where in the equation \eqref{eq3_9} the expectation has been taken with respect to the density $\xi_{\eta}(u_2,v_2|w^*) \propto v_2^{n_1+n_2-2} e^{-n_2 u_2 -
v_2 (1+w^*/\eta)}$. It can be easily seen that $\frac{\xi_{\eta}(u_2,v_2|w^*)}{\xi_{1}(u_2,v_2|w^*)}$ is non
increasing in $v_2$. Also we have $L^{\prime}(u_2-v_2c_{\eta}^*(w^*))$ is decreasing in $v_2$ as $c_{\eta}^*(w^*)>0$. $c_{\eta}(w)>0$ is guaranteed  as we have  $\phi_2(w^*)$ is a positive function. Consequently, using Lemma 3.4.2. in
\cite{lehmann2005testing}, we have
\[E_{\eta}\left(L^{\prime}\left(U_2 - V_2c^*_{1}(w^*)\right)\big|W^*=w^*\right) \ge
E_{1}\left(L^{\prime}\left(U_2 - V_2c^*_{1}(w^*)\right)\big|W^*=w^*\right)
= 0.\] The latter inequality means,
\begin{eqnarray}
c^*_{\eta}(w^*) \ge c^*_{1}(w^*),
\end{eqnarray}
because of \eqref{eq3_9}, which we have $E_{\eta}\left(L^{\prime}\left(U_2 - V_2c^*_{1}(w^*)\right)\big|W^*=w^*\right) \ge
E_{\eta}\left(L^{\prime}\left(U_2 - V_2c^*_{\eta}(w^*)\right)\big|W^*=w^*\right)$ and $E_{\eta}\left(L^{\prime}\left(U_2 - V_2c^*_{\eta}(w^*)\right)\big|W^*=w^*\right)$ is decreasing in $c^*_{\eta}(w^*)$.
Again $c^*_{1}(w^*)$ is the unique solution of
\begin{eqnarray}
\int_{0}^{\infty}\int_{0}^{\infty}L^{\prime}\left(u_2 - v_2c^*_{1}(w^*)\right)e^{-n_2u_2}v^{n_1+n_2-2}_2
e^{-v_2(1+w^*)}dv_2du_2=0.
\end{eqnarray}
By the transformation $u_2=u_2$, $z=v_2(1+w^*)$ we get
\begin{eqnarray}\label{sttheq34}
\int_{0}^{\infty}\int_{0}^{\infty}L^{\prime}\left(u_2 - \frac{z
c^*_1(w^*)}{1+w^*}\right) e^{-n_2 u_2}z^{n_1+n_2-2}e^{-z}dzdu_2=0.
\end{eqnarray}
Comparing (\ref{sttheq34}) with (\ref{bstar01}), we get
$c^*_1(w^*)=b_{02}(1+w^*)$. Moreover (\ref{c02}) and
(\ref{sttheqbstar01}) ensure that
\begin{eqnarray}
 b_{02}\ge c_{02}.
\end{eqnarray}
Consider a function
$\phi_{02}(w^*)=\max\left\{c_{02},b_{02}(1+w^*)\right\}$. Then we
have $c^*_{\eta}(w^*) > \phi_{02}(w^*)> c_{02}$ on a set of positive
probability. Since $R_2(c)$ is strictly bowl shaped, so for $c <
c^*_{\eta}(w^*)$, $R_2(c)$ is a decreasing function of $c$.
Therefore,
$R(\delta_{2S}) \le R(\delta_{02})$. This completes the proof of the theorem. \hfill\(\Box\)\\

\begin{corollary}\label{coro3_1}
Consider an estimator of the form (\ref{equiphi2}), that is
$\delta_{\phi_{2}}\left(\underline{X_{(1)}},\underline{T}\right)
=X_{2(1)} - T_{2}\phi_{2}\left(W^*\right).$ Define
$\phi_2^{*}(w^*)=\max\{\phi_{2}\left(w^*\right),b_{02}(1+w^*)\}$.
Then the estimator
$$\delta_{\phi_2^{*}}=X_{2(1)} - T_{2}\phi^{*}_{2}\left(W^*\right)$$
improve upon $\delta_{\phi_{2}}$ provided $P(\phi_{2}(W^*)<
b_{02}(1+W^*)) \ne 0$ under a bowl-shaped location invariant loss
function $L(t)$.
\end{corollary}

%
%
\begin{example}\rm
We consider the squared error loss function $L_1=\left(\frac{\delta-\mu_2}{\sigma_2}\right)^2$. Then the
estimator
$$\delta^Q_{2S}=X_{2(1)} - \max\left\{\frac{1}{n_2^2},\frac{1}{n_2(n_1+n_2-1)}(1+W^*)\right\}
T_2$$ improve upon the BAEE of $\mu_2$. It is emphasized that the improved estimator for $\mu_2$, $\delta^Q_{2S}$, is the same as the estimator discussed in Corollary 3.3 by \cite{vijayasree1995componentwise} for $k=2$.

The estimator which improve upon the UMVUE is obtained as
$$\delta^{Q}_{2MV}=X_{2(1)} - \max\left\{\frac{1}{n_2(n_2-1)},\frac{1}{n_2(n_1+n_2-1)}(1+W^*)\right\}
T_2$$.

The estimator $$ \delta^{Q}_{2ML}=X_{2(1)} - \frac{1}{n_2(n_1+n_2-1)}(1+W^*) T_2 $$ dominate the MLE of $\mu_2$.
\end{example}
\begin{example}
Consider the linex loss function $L_2=e^{a\left(\frac{\delta-\mu_2}{\sigma_2}\right)}-a\left(\frac{\delta-\mu_2}{\sigma_2}\right)-1$ and we have
for $n_2>a$,
$$b_{02}=\frac{1}{a}\left[\left(\frac{n_2}{n_2-a}\right)^{1/(n_1+n_2-1)} - 1\right]~~~\mbox{
and
}~~~c_{02}=\frac{1}{a}\left[\left(\frac{n_2}{n_2-a}\right)^{1/n_2} -1\right]$$ 
Then the  estimator
$$\delta^L_{2S}=X_{2(1)} - \min\left\{c_{02},b_{02}(1+W^*)\right\} T_2$$
improve upon the BAEE of $\mu_2$.  The estimator which improve upon
the UMVUE of $\mu_2$ is obtained as
$$\delta^{L}_{2MV}=X_{2(1)} - \max\left\{\frac{1}{n_2(n_2-1)},b_{02}(1+W^*)\right\} T_2$$ The estimator
$\delta^{L}_{2ML}=X_{2(1)} - b_{02}(1+W^*) T_2$
dominate the MLE of $\mu_2$.
\end{example}

In order to put more information in our estimating problem, we
consider another class of estimators of the form
\begin{eqnarray} \label{eq313}
\delta_{\phi_{2}^*}\left(\underline{X_{(1)}},\underline{T}\right)=X_{2(1)} - T_{2}
\phi_{2}^*\left(W^*,W_2\right),
\end{eqnarray}
where $W^*=\frac{T_1}{T_2}$, $W_2=\frac{X_{1(1)}}{T_2}$ and
$\phi_{2}^*$ is a measurable function.
\begin{theorem}
Let $b^*_{02}$ be the unique solution of
\begin{eqnarray}\label{bstar02}
E\left(L^{\prime}(U_2 - b^*_{02}Z_2)\right)=0,
\end{eqnarray}
where $Z_2$ follow a $Gamma(n_1+n_2,1)$ distribution and $U_2$
follow a standard exponential $Exp(1/n_2)$ distribution. Consider an
estimator of the form
\begin{eqnarray}
\delta^*_{2S}=X_{2(1)} -  \phi_{02}^*(W^*, W_2)T_2
\end{eqnarray}
with
\begin{eqnarray}
\phi^*_{02}(w^*)=\left\{\begin{array}{lcl}
\max\left\{c_{02},b^*_{02}(1+w^*+w_2)\right\} & , & w_2 > 0 \\
c_{02} & , & \text{otherwise}
\end{array}
\right.
\end{eqnarray}
where $c_{02}$ is given in (\ref{c02}). Then the risk of the
estimator $\delta^*_{2S}$ is nowhere larger than that of
$\delta_{02}$ under a bowl shaped invariant loss function $L(t)$
provided
\begin{eqnarray}\label{sttheq12}
E\left(L^{\prime}\left(U_2 - b^*_{02}V_2\right)V_2\right) \le 0.
\end{eqnarray}
\end{theorem}
\noindent \textbf{Proof:} The proof is analogous to Theorem
\ref{thmst31} and that's why is omitted.

\begin{corollary}\label{coro3_3}
Consider an estimator of the form (\ref{eq313}), that is
$\delta_{\phi_{2}}^*\left(\underline{X_{(1)}},\underline{T}\right)=X_{2(1)} - T_{2}\phi_{2}^*\left(W^*,W_2\right),$

Define $$\psi_2^{*}(w^*,w_2)=\left\{\begin{array}{lcl}

        \max\{\phi_{2}^*\left(w^*,w_2\right),b^*_{02}(1+w^*+w_2)\}, ~~~ w_2>0\\\\

        \phi_{2}^*\left(w^*,w_2\right),~~~~~~~~~~~~~~otherwise
    \end{array}
\right.$$

Then the estimator
        $$\delta_{\psi_2^{*}}=X_{2(1)} - T_{2}\psi^{*}_{2}\left(W^*,W_2\right)$$
        improve upon $\delta_{\phi_{2}}$ provided $P(\phi_{2}^*\left(w^*,w_2\right)< b^*_{02}(1+w^*+w_2)) \ne 0$
        under a bowl-shaped location invariant loss function $L(t)$.
    \end{corollary}
We can notice that the UMVUE and MLE of $\mu_2$ belong to the class
of estimators in (\ref{eq313}). As a consequence of Corollary
\ref{coro3_3} we get the following dominance result.

\begin{corollary}\label{corm3_4}
With respect to a bowl-shaped 
invariant loss function $L(t)$,
        the estimator $\delta_{2MV}$ is inadmissible and dominated by
        $$\delta^{I2}_{2MV}=\left\{\begin{array}{lcl}
        X_{2(1)} - \max\left\{\frac{1}{n_2(n_2-1)},b^*_{02}(1+W^*+W_2)\right\} T_2, W_2>0\\\\
        X_{2(1)}-\frac{1}{n_2(n_2-1)} T_2, ~~~~~~~~~~~~~~~~~~~~~~~~~otherwise
            \end{array}
    \right.,
        $$
        the estimator $\delta_{2ML}$ is inadmissible and dominated by $$\delta^{I2}_{2ML}=\left\{\begin{array}{lcl}
        X_{2(1)} - \max\left\{0,b^*_{02}(1+W^*+W_2)\right\} T_2,~~~ W_2>0\\\\
        X_{1(1)},~~~~~~~~~~~~~~~~~~~~~~~~otherwise
            \end{array}
    \right.,$$
\end{corollary} 
\begin{example}
Under the squared error loss function $L_1=\left(\frac{\delta-\mu_2}{\sigma_2}\right)^2$, the
    estimator of the form
    \begin{eqnarray}
        \delta_{2S}^{*,Q}=X_{2(1)} - \phi_{02}^*(W^*, W_2)T_2
    \end{eqnarray}
    with
    \begin{eqnarray}
        \phi_{02}^*(w^*,w_2)=\left\{\begin{array}{lcl}
            \max\left\{\frac{1}{n_2^2},\frac{(1+w^*+w_2)}{n_2(n_1+n_2)}\right\}, ~~ w_2 > 0 \\
            \frac{1}{n_2^2},~~~~~~~~~~~~~~~~~~~~~~~~~~~~\text{otherwise}
        \end{array}
        \right.,
    \end{eqnarray}
improves upon the BAEE of $\mu_2$.   The UMVUE of $\mu_2$ is
inadmissible and dominated by
$$\delta^{2,Q}_{2MV}=\left\{\begin{array}{lcl}
    X_{2(1)} - \max\left\{\frac{1}{n_2(n_2-1)},\frac{(1+W^*+W_2)}{n_2(n_1+n_2)}\right\} T_2, W_2>0\\\\
    X_{2(1)}-\frac{1}{n_2(n_2-1)} T_2, ~~~~~~~~~~~~~~~~~~~~~~~~~otherwise
\end{array}
\right.,
$$
Also, we have that $\delta_{2ML}$ is inadmissible and dominated by
$$\delta^{2,Q}_{2ML}=\left\{\begin{array}{lcl}
    X_{2(1)} - \frac{(1+w^*+w_2)}{n_2(n_1+n_2)} T_2,~~~ W_2>0\\\\
   X_{2(1)},~~~~~~~~~~~~~~~otherwis1e
\end{array}
\right.,$$
\end{example}
\begin{example}
Consider the linex loss function $L_2=e^{a\left(\frac{\delta-\mu_2}{\sigma_2}\right)}-a\left(\frac{\delta-\mu_2}{\sigma_2}\right)-1$. Under this
loss for $n_2>a$ function we have
 $$b^*_{02}=\frac{1}{a}\left[\left(\frac{n_2}{n_2-a}\right)^{1/(n_1+n_2)} - 1\right]~~~~\mbox{ and }
~~~c_{02}=\frac{1}{a}\left[\left(\frac{n_2}{n_2-a}\right)^{1/n_2} - 1\right]$$
The improved estimators for $\mu_2$ can be easily derived form
Corollary \ref{corm3_4} by substituting the values of $b^*_{02}$ and
$c_{02}$.
\end{example}

\subsection{A class of improved estimators for $\mu_2$}
In this section, following the technique of \cite{kubokawa1994unified}, we
derive a class of improved estimators. We consider a class of
estimator of the form
\begin{eqnarray}\label{2kubo}
\delta_{\varphi_2}=X_{2(1)} - T_2\varphi_2(W^*).
\end{eqnarray}

Let $U_2 = \dfrac{X_{2(1)} - \mu_2}{\sigma_2}$ with pdf $h(u_2) = n_2
e^{-n u_2}$, $u_2 > 0$ and $V_2 = \dfrac{T_2}{\sigma_2}$. The joint pdf of $V_2$, and $W^*$ is,
\begin{eqnarray*}
g_{\eta}(v_2,w^*) \propto \eta^{-n_1+1}e^{-v_2\left(1+\eta^{-1}
w^*\right)}(w^*)^{n_1-2}v_2^{n_1+n_2-3},~~v_2>0,w^*>0,0<\eta \le 1
\end{eqnarray*}
and $G_{\eta}(v_2, z) = \int_0^{z} g_{\eta}(v_2,x) dx$.

\begin{theorem}\label{thmkub2}
Let the function $\varphi_{2}(z)$ satisfies the following
assumptions
\begin{itemize}
\item[(i)] $\varphi_{2}(z)$ is non-increasing and
$\lim_{z\rightarrow\infty}\varphi_{2}(z)=c_{02}$,
 \item[(ii)]$\int_{0}^{\infty}
\int_{0}^{\infty}L^{\prime}\left(u_2 - v_2\varphi_{2}(z)\right)h(u_2)v_2G_{1}(v_2,z)du_2
dv_2 \le 0.$
\end{itemize}
Then the risk of the estimator $\delta_{\varphi_{2}}$ given in
(\ref{2kubo}) dominates $\delta_{02}$ under a bowl shaped location-scale
invariant loss function $L(t)$.
\end{theorem}

\noindent{\textbf{Proof}:}  Proof of this result is similar to Theorem \ref{thmkubo1}.

	\begin{remark} \label{rembz2}
		The limit case in Theorem \ref{thmkub2} provides the \cite{brewster1974improving} - type estimator (see \cite{kubokawa1994unified} for details). 
		That means, $\delta_{\varphi_{2BZ}}=X_{2(1)}-T_2\varphi_{BZ}(W^*)$ is the \cite{brewster1974improving}- type estimator of $\mu_2$, where $\varphi_{2BZ}$ is given in Equation \eqref{eqbz2}
		\begin{equation} \label{eqbz2}
			\int_{0}^{\infty} \int_{0}^{\infty}L^{\prime}\left(u_2 - v_2\varphi_{2BZ}(z)\right) h(u_2) v_2 G_{1}(v_2,z) du_2 dv_2 = 0.
		\end{equation}
		As in Remark \ref{rembz1}, it can be proved that this estimator is a generalized Bayes estimator for $\mu_2$ under the improper prior $\pi(\mu_2, \sigma_1, \sigma_2) = \frac{1}{\sigma_1 \sigma_2}$, $\mu_2 \in R$, $0< \sigma_1 < \sigma_2$ and the bowl-shaped loss function $L(t)$.
	\end{remark}

 \begin{example} Suppose the function $\varphi_{2}(z)$
satisfies the following assumptions
\begin{itemize}
    \item[(i)] $\varphi_{2}(z)$ is non-increasing and $\lim_{z\rightarrow\infty}\varphi_{2}(z)= \frac{1}{n_2^2}$,
    \item[(ii)]$\varphi_{2}(z) \ge\varphi_{2BZ}(z)$,
    with
    $$\varphi_{2BZ}(z)
   	=\frac{1}{n_2(n_1+n_2-1)}\frac{\int_{0}^{z}\frac{x^{n_1-2}}{(1+x)^{n_1+n_2-1}}}{\int_{0}^{z}\frac{x^{n_1-2}}{(1+x)^{n_1+n_2}}}$$
\end{itemize}
Then the risk of the estimator $\delta_{\varphi_{2}}$ given in
(\ref{2kubo}) dominates $\delta_{02}$ under squared error loss
$L_1\left(\frac{\delta-\mu_2}{\sigma_2}\right)$.
\end{example}

\begin{example} Suppose the function $\varphi_{2}(z)$
satisfies the following assumptions
\begin{itemize}
    \item[(i)] $\varphi_{2}(z)$ is non-decreasing and
 $\lim_{z\rightarrow\infty}\varphi_{2}(z)= \frac{1}{a}\left[\left(\frac{n_2}{n_2-a}\right)^{1/n_2} - 1\right]$,
    \item[(ii)] $\int_{0}^{\infty} \int_{0}^{\infty} a \left(e^{a (u_2 - v_2\varphi_{2}(z) )} -1\right)h(u_2)v_2G_{2}(v_2,z)du_2 dv_2  \le 0.$
\end{itemize}
Then the risk of the estimator $\delta_{\varphi_{2}}$ given in
(\ref{2kubo}) dominates $\delta_{02}$ under linex loss $L_2\left(\frac{\delta-\mu_2}{\sigma_2}\right)$.
\end{example}




\section{Improved estimation under Generalized Pitman Closeness}\label{sec4}

In this section we consider the estimation problem under the generalized Pitman nearness (GPN) criterion.  For a discussion Pitman nearness criterion we refer to \cite{garg2023unified}.The generalized Pitman nearness (GPN) criterion (see \cite{nayak1990estimation} and \cite{kubokawa1991equivariant}) based on general loss function $W(\underline{\theta},\delta)$ is defined as follows:
\begin{definition}
	Let $\underline{X}$ be a random vector having a distribution involving an unknown parameter $\underline{\theta}\in \Theta$ and let $\tau(\underline{\theta})$ be a real-valued estimand. Let $\delta_1$ and $\delta_2$ be two estimators of the estimand $\tau(\underline{\theta})$. Also, let $W(\underline{\theta},a)$ be a specified loss function for estimating $\tau(\underline{\theta})$. Then, the generalized Pitman nearness (GPN) of $\delta_1$ relative to $\delta_2$ is defined by
	$$GPN(\delta_1,\delta_2;\underline{\theta})=P_{\underline{\theta}}[W(\underline{\theta},\delta_1) <W(\underline{\theta},\delta_2)]+\frac{1}{2} P_{\underline{\theta}}[W(\underline{\theta},\delta_1) =W(\underline{\theta},\delta_2)], \; \; \underline{\theta}\in\Theta.$$
	The estimator $\delta_1$ is said to be nearer to $\tau(\underline{\theta})$ than $\delta_2$, under the GPN criterion, if $GPN(\delta_1,\delta_2;\underline{\theta})\geq \frac{1}{2},\;\forall\; \underline{\theta}\in\Theta$, with strict inequality for some $\underline{\theta}\in\Theta$.
\end{definition}
The following lemma, taken from \cite{garg2023unified}, will be useful in proving the main results of this section (\cite{nayak1990estimation} and \cite{zhou2012pitman}).
\begin{lemma}
Let $Y$ be a r.v. having the Lebesgue p.d.f., let $X$ be a positive r.v. having the Lebesgue p.d.f. and let $m_{YX}$ be the median of r.v. $\frac{Y}{X}$. Let $L:\Re\rightarrow [0,\infty)$ be a function such that $L(t)$ is strictly decreasing on $(-\infty,0)$ and strictly increasing on $(0,\infty)$. Then, for $-\infty< c_1<c_2\leq m_{YX}$ or $-\infty<m_{YX}\leq c_2<c_1$,
$GPN= P[L(Y-c_2X)<L(Y-c_1X)]+\frac{1}{2} P[L(Y-c_2X)=L(Y-c_1X)]>\frac{1}{2}$.
\end{lemma}

Note that, any affine equivariant estimator is of the form $\delta_{i,c}(\underline{X})=X_{i(1)}-c T_i,\; c>0,\; i=1,2$. An immediate consequence of Lemma 4.1 is that the Pitman nearest affine equivariant estimator (PNAEE) of $\mu_i$, under the GPN criterion, is $\delta_{i,P}(\underline{X})=X_{i(1)}-m_{0,i}T_i$, where $m_{0,i}$ is the median of the r.v. $\frac{U_i}{V_i}=\frac{X_{i(1)}-\mu_i}{T_i},\;i=1,2$. Note that, in this section we use $\underline{\theta}=(\mu_{1},\mu_2,\sigma_{1},\sigma_{2})$ and
$\Theta_0=\{(\mu_{1},\mu_2,\sigma_{1},\sigma_{2}): \mu_{1}\in \Re, \mu_2\in \Re, 0<\sigma_{1}\leq \sigma_{2}<\infty\}$.
\vspace*{2mm}

\subsection{\textbf{Improved estimation of $\mu_{1}$}} 

Let $\delta_{\xi}(\underline{X})=X_{1(1)}-T_1 \xi(W)$ and $\delta_{\psi}(\underline{X})=X_{1(1)}-T_1\psi(W)$ be two equivariant estimators of $\mu_1$ (see \eqref{equiphi1}), where $\xi$ and $\psi$ are real-valued positive functions. Define $U_1=\frac{X_{1(1)}-\mu_1}{\sigma_{1}}$, $V_1=\frac{T_1}{\sigma_1}$ and $\eta=\frac{\sigma_1}{\sigma_2}\leq 1$. Then, the GPN of $\delta_{\xi}(\underline{X})$ relative to $\delta_{\psi}(\underline{X})$ is given by 
\begin{align*}
	GPN(\delta_{\xi},\delta_{\psi};\underline{\theta})&=P[L(U_1-V_1\xi(W)) <L(U_1-V_1\psi(W))]
 +\frac{1}{2} P[L(U_1-V_1\xi(W)) =L(U_1-V_1\psi(W))]\\
	&=E^{W}[P[L(U_1-V_1\xi(W)) <L(U_1-V_1\psi(W))\vert W]] \\
	 &+\frac{1}{2} E^W[P[L(U_1-V_1\xi(W)) =L(U_1-V_1\psi(W))]\vert W].
\end{align*}
For fixed $w>0$, define 
\begin{align} \label{eq:P1}
	R_{1,\eta}(\xi(w),\psi(w),w)&= P_{\underline{\theta}}[L(U_1-V_1\xi(w)) <L(U_1-V_1\psi(w))\vert W=w] \nonumber
	\\&\qquad +\frac{1}{2}P_{\underline{\theta}}[L(U_1-V_1\xi(w))\! =\!L(U_1-V_1\psi(w))\vert W=w].
\end{align}
For any fixed $0<\eta\leq 1$ and $w>0$, the conditional distribution of $\frac{U_1}{V_1}$ given $W=w$ is
\begin{equation*}
	g(x)= \begin{cases}
	\frac{n_1(n_1+n_2-2)(n_1x+1+\eta \,w)^{-n_1-n_2+1}}{(1+\eta \,w)^{-n_1-n_2+2}},&\text{if}\;\; x>0\\
	0, &\text{otherwise}
\end{cases},
\end{equation*}
and the median of the conditional distribution of $\frac{U_1}{V_1}$ given $W=w$ is
\begin{equation}
	m_{\eta}(w)= \frac{1}{n_1}\left(2^{\frac{1}{n_1+n_2-2}}-1\right)(1+\eta w).
\end{equation}
Since, for every fixed $w>0$, $m_{\eta}(w)$ is an increasing function of $\eta$. Therefore, for $w>0$, we have
\begin{align}
\label{eq:P2}	l(w)& =\inf_{\eta\leq 1} m_{\eta}(w) =\lim_{\eta\to 0} m_{\eta}(w)=\frac{1}{n_1}(2^{\frac{1}{n_1+n_2-2}}-1)\\
	\;\text{and}\;\;
	u(w)&=\sup_{\eta\leq 1} m_{\eta}(w)=\lim_{\eta\to 1}
\label{eq:P3}	 m_{\eta}(w)=\frac{1}{n_1}(2^{\frac{1}{n_1+n_2-2}}-1)(1+w)
\end{align}

\vspace*{2mm}

Now we develop the following theorem that provides an estimator dominating an arbitrary equivariant estimator $\delta_{\psi}(\underline{X})=X_{1(1)}-T_1\psi(W)$.
\vspace*{2mm}

\begin{theorem}\label{pitth1}
	Let $\delta_{\psi}(\underline{X})=X_{1(1)}-T_1\psi(W)$ be a equivariant estimator of $\mu_1$. Let $l(w)$ and $u(w)$ be functions as defined by \eqref{eq:P2} and \eqref{eq:P3}, respectively, such that $l(w)\leq m_{\eta}(w)\leq u(w), \;\forall\;\eta\leq 1$ and any $w>0$. For any fixed $w>0$, define $\psi^{*}(w)\!=\!\max\{l(w),\min\{\psi(w),u(w)\}\}$. Then, under the GPN criterion with a general loss function,
the estimator $\delta_{\psi^{*}}(\underline{X})\!=\!X_{1(1)}-T_1\psi^{*}(W)$ is Pitman nearer to $\mu_1$ than the estimator $\delta_{\psi}(\underline{X})=X_{1(1)}-T_1\psi(W)$, for all $\underline{\theta}\in \Theta_0$, provided $P[l(W)\leq \psi(W)\leq u(W)]<1,\;\forall\; \underline{\theta}\in \Theta_0$.
\end{theorem}
\noindent{\textbf{Proof}:}
	The GPN of the estimator $\delta_{\psi^{*}}(\underline{X})=X_{1(1)}-T_1\psi^{*}(W)$ relative to $\delta_{\psi}(\underline{X})=X_{1(1)}-T_1\psi(W)$ can be written as 
	\begin{align*}
		GPN(\delta_{\psi^{*}},\delta_{\psi};\underline{\theta})=E^W[R_{1,\eta}(\xi(W),\psi(W),W)],\;\; \eta\leq 1,
	\end{align*}
	where $R_{1,\eta}(\cdot,\cdot,\cdot)$ is defined by \eqref{eq:P1}.
	\vspace*{2mm}
	
	Let $A=\{w:\psi(w)<l(w)\}$, $B=\{w:l(w)\leq \psi(w)\leq u(w)\}$ and $C=\{w:\psi(w)>u(w)\}$. Clearly
	$$\psi^{*}(w)=\begin{cases} l(w), & w\in A\\  \psi(w), & w\in B\\ u(w), & w\in C\end{cases}.$$
	Since $l(w)\leq m_{\eta}(w)\leq u(w),\; \forall\; \eta\leq 1$ and $w>0$, using Lemma 4.1, we have $R_{1,\eta}(\psi^{*}(w),\psi(w),w)> \frac{1}{2}, \; \forall \; \eta\leq 1$, provided $w\in A\cup C$. Also, for $w\in B$, $R_{1,\eta}(\psi^{*}(w),\psi(w),w)= \frac{1}{2}, \; \forall \; \eta\leq 1$. Let $f_W(w\vert \eta)$ be the p.d.f. of r.v. $W$. Since $P[A\cup C]>0,\; \forall \; \underline{\theta}\in \Theta_0$, we have
	\begin{align*}
GPN(\delta_{\psi^{*}},\delta_{\psi};\underline{\theta})
		&\!=\!\int_{A}\! R_{1,\eta}(\psi^{*}(w),\psi(w),w) f_W(w\vert \eta) dw \!+ \!\int_{B}\! R_{1,\eta}(\psi^{*}(w),\psi(t),w) f_W(w\vert \eta) dw \!\\&\qquad\qquad\qquad\qquad+\! \int_{C}\! R_{1,\eta}(\psi^{*}(w),\psi(w),w) f_W(w\vert \eta) dw\\
		&> \frac{1}{2}, \;\; \underline{\theta}\in\Theta_0.
	\end{align*}

	The following corollary provides improvements over the unrestricted Pitman nearness affine equivariant estimator (PNAEE) $\delta_{1,P}(\underline{X})=X_{1(1)}- m_{0,1}T_1$, where $m_{0,1}=\frac{1}{n_1}(2^{\frac{1}{n_1-1}}-1)$ is the median of r.v. $\frac{U_1}{V_1}$.
 \begin{corollary}
 	Let $l(w)$ and $u(w)$ be as defined by \eqref{eq:P2} and \eqref{eq:P3}, respectively. Suppose that $P[l(W)\leq m_{0,1}\leq u(W)]<1,\; \forall\; \underline{\theta}\in \Theta_0$. For any fixed $w>0$, define  $\xi^{*}(w)\!=\!\max\{l(w),\min\{m_{0,1},u(w)\}\}=\min\{m_{0,1},u(w)\}.$ Then, for every $\underline{\theta}\in \Theta_0$,
	the estimator $\delta_{\xi^{*}}(\underline{X})\!=\!X_{1(1)}- \xi^{*}(W)T_1$ is Pitman nearer to $\mu_1$ than the PNAEE $\delta_{1,P}(\underline{X})=X_{1(1)}- m_{0,1}T_1$, under the GPN criterion.
 \end{corollary}

		The following corollary provides improvements over the unrestricted Best affine equivariant estimator (BAEE) $\delta_{1,P}(\underline{X})=X_{1(1)}- c_{01}T_1$, where $c_{0,1}$ is as defined by \eqref{c01}.
\begin{corollary}
	Let $l(w)$ and $u(w)$ be as defined by \eqref{eq:P2} and \eqref{eq:P3}, respectively. Suppose that $P[l(W)\leq c_{01}\leq u(W)]<1,\; \forall\; \underline{\theta}\in \Theta_0$. Define, for any fixed $w$, $\xi^{*}(w)\!=\!\max\{l(w),\min\{c_{01},u(w)\}\}$. Then, for every $\underline{\theta}\in \Theta_0$,
	the estimator $\delta_{\xi^{*}}(\underline{X})\!=\!X_{1(1)}-T_1 \xi^{*}(W)$ is Pitman nearer to $\mu_1$ than the BAEE $\delta_{01}(\underline{X})=X_{1(1)}- c_{01}T_1$, under the GPN criterion.
\end{corollary}

	\subsection{\textbf{Improved estimation of $\mu_{2}$}} 
	
	Let $\delta_{\xi}(\underline{X})=X_{2(1)}-T_2 \xi(W^{*})$ and $\delta_{\psi}(\underline{X})=X_{2(1)}-T_2\psi(W^{*})$ be two equivariant estimators of $\mu_2$ (see \eqref{equiphi2}), where $\xi$ and $\psi$ are real-valued positive functions. Define $U_2=\frac{X_{2(1)}-\mu_2}{\sigma_{2}}$, $V_2=\frac{T_2}{\sigma_2}$ and $\eta=\frac{\sigma_1}{\sigma_2}\leq 1$. Then, the GPN of $\delta_{\xi}(\underline{X})$ relative to $\delta_{\psi}(\underline{X})$ is given by 
	\begin{align*}
		GPN(\delta_{\xi},\delta_{\psi};\underline{\theta})&=P[L(U_2-V_2\xi(W^{*})) <L(U_2-V_2\psi(W^{*}))]
		\\&\quad +\frac{1}{2} P[L(U_2-V_2\xi(W^{*})) =L(U_2-V_2\psi(W^{*}))]\\
		&=E^{W^{*}}[P[L(U_2-V_2\xi(W^{*})) <L(U_2-V_2\psi(W^{*}))\vert W^{*}]]
		\\&\quad +\frac{1}{2} E^{W^{*}}[P[L(U_2-V_2\xi(W^{*})) =L(U_2-V_2\psi(W^{*}))]\vert W^{*}].
	\end{align*}
	For fixed $w^{*}>0$, define 
	\begin{align} \label{eq:P31}
		R_{2,\eta}(\xi(w^{*}),\psi(w^{*}),w^{*})&= P_{\underline{\theta}}[L(U_2-V_2\xi(w^{*})) <L(U_2-V_2\psi(w^{*}))\vert W^{*}=w^{*}] \nonumber
		\\&\qquad +\frac{1}{2}P_{\underline{\theta}}[L(U_2-V_2\xi(w^{*}))\! =\!L(U_2-V_2\psi(w^{*}))\vert W^{*}=w^{*}].\end{align}

	For any fixed $0<\eta\leq 1$ and $w^{*}>0$, the conditional distribution of $\frac{U_2}{V_2}$ given $W^{*}=w^{*}$ is
	\begin{equation*}
		g(x)= \begin{cases}
			\frac{n_2(n_1+n_2-2)\left(n_2x+1+\frac{w^{*}}{\eta}\right)^{-n_1-n_2+1}}{\left(1+\frac{w^{*}}{\eta}\right)^{-n_1-n_2+2}},&\text{if}\;\; x>0\\
			0, &\text{otherwise}
		\end{cases},
	\end{equation*}
	and the median of the conditional distribution of $\frac{U_2}{V_2}$ given $W^{*}=w^{*}$ is
	\begin{equation*}
		m_{\eta}(w^{*})= \frac{1}{n_2}\left(2^{\frac{1}{n_1+n_2-2}}-1\right)\left(1+\frac{ w^{*}}{\eta}\right).
	\end{equation*}
	Since, for every fixed $w^{*}>0$, $m_{\eta}(w^{*})$ is a decreasing function of $\eta$. Therefore, for $w^{*}>0$, we have
	\begin{align}
		\label{eq:P4}	l(w^{*})& =\inf_{\eta\leq 1} m_{\eta}(w^{*}) =\lim_{\eta\to 1} m_{\eta}(w^{*})=\frac{1}{n_2}\left(2^{\frac{1}{n_1+n_2-2}}-1\right)(1+w^{*})\\
		\;\text{and}\;\;
		u(w^{*})&=\sup_{\eta\leq 1} m_{\eta}(w^{*})=\lim_{\eta\to 0}
		\label{eq:P5}	 m_{\eta}(w^{*})=\infty
	\end{align}
	
Now we develop the following theorem that provides an estimator dominating an arbitrary equivariant estimator $\delta_{\psi}(\underline{X})=X_{2(1)}-T_2\psi(W^{*})$.
\begin{theorem}
Let $\delta_{\psi}(\underline{X})=X_{2(1)}-T_2\psi(W^{*})$ be a equivariant estimator of $\mu_2$. Let $l(w^{*})$ and $u(w^{*})$ be functions as defined by \eqref{eq:P4} and \eqref{eq:P5}, respectively, such that $l(w^{*})\leq m_{\eta}(w^{*})\leq u(w^{*}), \;\forall\;\eta\leq 1$ and any $w^{*}>0$. For any fixed $w^{*}>0$, define $\psi^{*}(w^{*})\!=\!\max\{l(w^{*}),\min\{\psi(w^{*}),u(w^{*})\}\}$. Then, under the GPN criterion with a general loss function,
	the estimator $\delta_{\psi^{*}}(\underline{X})\!=\!X_{2(1)}-T_2\psi^{*}(W^{*})$ is Pitman nearer to $\mu_2$ than the estimator $\delta_{\psi}(\underline{X})=X_{2(1)}-T_2\psi(W^{*})$, for all $\underline{\theta}\in \Theta_0$, provided $P[l(W^{*})\leq \psi(W^{*})\leq u(W^{*})]<1,\;\forall\; \underline{\theta}\in \Theta_0$.
\end{theorem}
\noindent{\textbf{Proof}:} Proof of this theorem is similar to the proof of the Theorem \ref{pitth1}.

\vspace{0.5cm}
The following corollary provides improvements over the unrestricted Pitman nearness affine equivariant estimator (PNAEE) $\delta_{2,P}(\underline{X})=X_{2(1)}- m_{0,2}T_2$, where $m_{0,2}=\frac{1}{n_2}(2^{\frac{1}{n_2-1}}-1)$ is the median of r.v. $\frac{U_2}{V_2}$.
	\begin{corollary}
 Let $l(w^{*})$ and $u(w^{*})$ be as defined by \eqref{eq:P4} and \eqref{eq:P5}, respectively. Suppose that $P[l(W^{*})\leq m_{0,2}\leq u(W^{*})]<1,\; \forall\; \underline{\theta}\in \Theta_0$.  For any fixed $w^{*}>0$, define  $\xi^{*}(w^{*})\!=\!\max\{l(w^{*}),m_{0,2}\}.$ Then, for every $\underline{\theta}\in \Theta_0$,
	the estimator $\delta_{\xi^{*}}(\underline{X})\!=\!X_{2(1)}- \xi^{*}(W^{*})T_2$ is Pitman nearer to $\mu_2$ than the PNAEE $\delta_{2,P}(\underline{X})=X_{2(1)}- m_{0,2}T_2$, under the GPN criterion.
	\end{corollary}
The following corollary provides improvements over the unrestricted BAEE $\delta_{2,P}(\underline{X})=X_{2(1)}- c_{02}T_2$, where $c_{02}$ is as defined by \eqref{c02}.
	\begin{corollary}
		Let $l(w^{*})$ and $u(w^{*})$ be as defined by \eqref{eq:P4} and \eqref{eq:P5}, respectively. Suppose that $P[l(W^{*})\leq c_{02}\leq u(W^{*})]<1,\; \forall\; \underline{\theta}\in \Theta_0$. Define, for any fixed $w^{*}$, $\xi^{*}(w^{*})\!=\!\max\{l(w^{*}),c_{02}\}$. Then, for every $\underline{\theta}\in \Theta_0$,
	the estimator $\delta_{\xi^{*}}(\underline{X})\!=\!X_{2(1)}-T_2 \xi^{*}(W^{*})$ is Pitman nearer to $\mu_2$ than the BAEE $\delta_{02}(\underline{X})=X_{2(1)}- c_{02}T_2$, under the GPN criterion.
	\end{corollary}


\section{Application to special sampling schemes}\label{sec5}
\subsection{Type II censoring}

In the realm of reliability and life-testing experiments, various situations arise where test units are either lost or withdrawn from the experiment before they fail. For instance, in industrial trials, units may break unexpectedly prior to the intended time, and in clinical trials, participants might drop out or the experiment could halt due to budget constraints. In many cases, the decision to remove units before they fail is premeditated, primarily to save time and reduce testing expenses. The data collected from such experiments are referred to as censored data, with one common type known as Type-II censoring. Under this scheme, the experimenter decides to conclude the experiment after a specified number of items ($r$, where $r$ is less than or equal to $n$) experience failure. For more detailed information on this subject, we  consulted reference \cite{balakrishnan2000progressive}.

Consider drawing a sample of size $n$ from an exponential distribution $E(\mu_i, \sigma _i)$, where the observations are arranged in ascending order, i.e., $X_{i(1)}^c \le X_{i(2)}^c \le \dots \le X_{i(n)}^c$ for ${ i = 1,2}$. Here, $X_{i(j)}^c$ represents the $j$th smallest observation in a sample of $n_i$ observations taken from an exponential $E(\mu_i, \sigma_i )$ population. Now, let's focus on the first $r$ ordered observations, which are $X_{i(1)}^c, X_{i(2)}^c,\dots, X_{i(r)}^c$, where $r$ is less than or equal to $n_i$, and $i = 1,2$. We aim to estimate $\mu_i$, based on type-II censored sample using the affine-invariant loss function $L(t)$, while adhering to the order constraint that $\sigma_1 \le \sigma_2$.\\
Define $T_i$ as $T_{i}=\sum_{j=1}^{n_{i}}\left(X_{ij}^c-X_{i(1)}^c\right)$.
In this setup, $(X_{1(1)}^c, X_{2(1)}^c, T_1,T_2)$ constitutes a minimal sufficient statistic. So we get the model  (\ref{model1}). Consequently, improved estimators of $\mu_{i}$ can be derived from the Theorems proved in Section \ref{sec2s} and \ref{sec3}
\subsection{Progressive type II censoring}
Lifetime distributions analyzed under censored sampling methods are widely utilized in various fields, including science, engineering, social sciences, public health, and medicine due to their broad applicability. Censored data falls into the category of progressive type-II censoring when they are managed by removing a predetermined number of surviving units at the time of failure of an individual unit.  This can be described as follows. In a life-testing experiment, $n$ units are involved, but only $m$ (where $m < n$) are observed continuously until they fail. Censoring takes place progressively in $m$ stages.  These $m$ stages provide failure times for the $m$ units that are entirely observed. At the moment of the first failure (the first stage), $R_1$ out of the $n-1$ remaining units are randomly withdrawn from the experiment. Subsequently, at the time of the second failure (the second stage), $R_2$ out of the $n-2-R_1$ remaining units are randomly withdrawn. Finally, when the $m$-th failure occurs (the $m$-th stage), all the remaining $R_m = n - m - R_1 - R_2 - \dots - R_m-1$ units are withdrawn. for more details see \cite{madi2010note},\cite{balakrishnan2007progressive} and \cite{balakrishnan2000progressive}.

In a life-testing experiment, there are $m_i$ independent units, each with a lifetime denoted as $X_1^p, X_2^p,\dots, X_{m_{i}}^p$. These lifetimes follow an exponential $E(\mu_{i},\sigma_{i})$. As the experiment progresses, when the $kth$ unit fails (for $k = 1, 2,\dots,n_i$ where $n_i\le m_i)$, a predetermined number of $R_k$ surviving units are taken out of the experiment. This setup allows for the systematic removal of surviving units at each failure event in the experiment. Where $R_n = m-n-R_1-R_2 -\dots-R_{n-1}$. Let we have two sets of data: one set is denoted as $X_{1:n_i:m_i}^p, X_{2:n_i:m_i}^p,\dots, X_{n_i:n_i:m_i}^p$, representing the progressive type-II censored sample of lifetimes for $m_i$ independent units, following an exponential distribution $E(\mu_{1},\sigma_{1})$. The other set is denoted as $Y_{1,n_i:m_i}^p, Y_{2,n_i:m_i}^p,\dots, Y_{n_i:n_i:m_i}^p$, representing the lifetimes of another set of $m_i$ independent units following an exponential distribution $E(\mu_2 ,\sigma_{2})$ associated restriction $\sigma_{1}\le \sigma_{2}.$

Our goal is to estimate the parameter $\mu_i$ based on these progressive type-II censored samples. To estimate this, we define  here $(X_{1,n_i:m_i}^p,T_1)$ as complete suffiecient statistic for and $(Y_{1,n_i:m_i}^p,T_2)$ are complete sufficient statistics. Where $T_1 = \sum_{j=1}^{m_i}(R_i+1)X_{j,n_i:m_i}^p-X_{1,n_i:m_i}^p$ and $T_2= \sum_{j=1}^{m_i}(R_i+1)Y_{j,n_i:m_i}^p-Y_{1,n_i:m_i}^p)$.
It's then shown that $X_{1:n:m}^p$ follows $E(\mu_1,\frac{\sigma_1}{m_i})$, $Y_{1:m:n}^p$ follows $E(\mu_2,\frac{\sigma_2}{m_i})$ for $i = 1, 2$ and $T_1$, $T _2$ follows a gamma distribution $\Gamma(m_i-1, \sigma_1)$, $\Gamma(m_i-1,\sigma_{2})$ respectively.

This setup leads to model (\ref{model1}), which allows us to derive improved estimators for $\mu_{i}$ based on the information provided in theorems proved in Section \ref{sec2} and \ref{sec3}. These theorems likely contain specific formulas or procedures for estimating $\mu_{i}$ in this context, which can be applied to our data.
%

\subsection{Record values}
Record values have a wide range of applications across diverse fields like hydrology, meteorology, and stock market analysis. Extensive research has been conducted in the scientific literature regarding records of independent and identically distributed (i.i.d.) random variables, along with their associated properties.
Consider a sequence of random variables $V_1, V_2,\dots, V_n$ that are both independent and identically distributed (iid). These random variables have a cumulative distribution function (CDF) denoted as $F(x)$and a probability density function (PDF) denoted as $f(x)$. We define a new sequence of random variables such that  for $k = 1$ set  $U(1) = 1$ and $U(k)$=min$\{j|j>U(k-1),V_j>V_{u(k-1)}\}$. This new sequence $\{X_k = V_{U(k)}, k \ge 1\}$, represents a sequence of maximal record statistics. For more discussion on this one may refer to \cite{madi2008improved} and \cite{bobotas2011improved}.

Suppose having two samples of record values: $X_1^r, X_2,^r\dots, X_{n}^r$ drawn from an $E(\mu_1,\sigma_{1})$population and $Y_1^r, Y_2^r,\dots, Y_{m}^r$ drawn from another exponential population $E(\mu_{2},\sigma_{2})$.
In the context of these sample sets, we are interested in estimating $\mu_i$, with respect to a general affine-invariant loss function $L(t)$, subject to the constraint that $\sigma_{1} \le \sigma_2$. Here, $(X_1^r,Y_1^r,T_1,T_2)$ are sufficient statistics where $T_{1}=\left(X_{(n)}^r-X_{(1)}^r \right)$  and follows $\Gamma(n-1,\sigma_{1})$, $T_{2}=\left(Y_{(m)}^r-Y_{(1)}^r\right)$ and follows $\Gamma(m-1,\sigma_{2})$. It's important to note that $X_{1}^r$ and $Y_{1}^r$ follows  exponential distribution $E(\mu_{1},\sigma_{1})$ and $E(\mu_{2},\sigma_{2})$ respectively. 
Therefore, we have a model that is similar to the one presented in (\ref{model1}), and the results discussed in Section \ref{sec2} are suited  here as well. On observing the model \ref{model1} we can state that improvement results parallel to section \ref{sec2s}, \ref{sec3} and \ref{sec4} can be obtained for these cases too.
\section{Simulation Study}\label{sec6}
In this section, we have conducted a simulation study to compare the risk performance of the improved estimators with respect to the squared error loss function. We have generated 20000 random samples of size $n_1$ and $n_2$ form $Exp(\mu_1,\sigma_1)$ and  $Exp(\mu_2,\sigma_2)$ respectively.  To compare the risk performance of the improved estimators of $\mu_1$, we calculate percentage risk improvement (PRI) with respect to   BAEE and UMVUE for different values of $\mu_1$, $\mu_2$, $\sigma_1$, $\sigma_2$ and $(n_1,n_2)$. The PRI of an estimator $\delta$  with respect to $\delta_0$ is defined as
\begin{eqnarray*}
	PRI(\delta)=\frac{Risk(\delta_0)-Risk(T)}{Risk(\delta_0)} \times 100
\end{eqnarray*}
In Table \ref{tab1_mu_1l1} and \ref{tab1_mu_1l2} we have tabulated the PRI with respect to BAEE, and in Table \ref{tab1_mu_1l2} and \ref{tab1_mu_1l4} we have tabulated the PRI with respect to UMVUE. From the tabulated value, we have the following observations. 
\begin{itemize}
	\item [(i)]  From the simulated value, we have seen that risk performance of $\delta_{1s}^Q$,  $\delta_{1s}^{*Q}$, $\delta_{1BZ}$ always better than BAEE.
	\item [(ii)] The PRI values of $\delta_{1s}^Q$,  $\delta_{1s}^{*Q}$ decreases as $\frac{\sigma_1}{\sigma_2}$ decreases. 
	\item [(iii)] The PRI values of $\delta_{1BZ}$ increases and then decreases. 
	\item[(iv)] From the simulated PRI, we have seen that sometimes BAEE performs better than  $\delta_{1MV}^{Q}$ and $\delta_{1MV}^{2Q}$ and sometimes  $\delta_{1MV}^{Q}$ and $\delta_{1MV}^{2Q}$ perform better than BAEE. So  there is no clear comparison between the estimators  $\delta_{1MV}^{Q}$, $\delta_{1MV}^{2Q}$ with BAEE. 
	\item[(v)] The estimators $\delta_{1MV}^{Q}$ and $\delta_{1MV}^{2Q}$ improves upon UMVUE. The PRI of  $\delta_{1MV}^{Q}$ and $\delta_{1MV}^{2Q}$ decreases as $\frac{\sigma_{1}}{\sigma_2}$ decreases.  
\end{itemize}
For the improved estimators of $\mu_2$, we have similar observations. 
\pagebreak
\begin{table}[h!]
		\centering
		\caption{Percentage risk improvement of the improved estimators of $\mu_1$ over the BAEE under quadratic loss function $L_1(t)$}
		\vspace{0.2cm}
		\label{tab1_mu_1l1}
\begin{tabular}{ccccccc|ccccc}
	\hline \hline\vspace{0.3cm}  
	$\sigma_1$           & $\sigma_2$ & \multicolumn{1}{c}{$\delta_{1s}^Q$} & \multicolumn{1}{c}{$\delta_{1S}^{*Q}$} & \multicolumn{1}{c}{$\delta_{1MV}^Q$} & \multicolumn{1}{c}{$\delta_{1MV}^{2Q}$} & $\delta_{1BZ}$ & \multicolumn{1}{c}{$\delta_{1s}^Q$} & \multicolumn{1}{c}{$\delta_{1S}^{*Q}$} &  \multicolumn{1}{c}{$\delta_{1MV}^Q$} & \multicolumn{1}{c}{$\delta_{1MV}^{2Q}$} & $\delta_{1BZ}$  \\\hline\hline \vspace{0.3cm}  
	 &        &   \multicolumn{5}{c|}{$(n_1,n_2)=(4,5),   ~~(\mu_1,\mu_2)=(0.1,0.3)$}                                                                                                                                                                                                     & \multicolumn{5}{c}{$(n_1,n_2)=(6,11), (\mu_1,\mu_2)=(0.1,0.3)$}                                                                                                                                                                                                  \\ \hline\hline
	\multirow{3}{*}{0.4} & 0.6        & \multicolumn{1}{c|}{2.14}            & \multicolumn{1}{c|}{2.50}               &  \multicolumn{1}{c|}{4.31}             & \multicolumn{1}{c|}{5.27}                & 1.46           & \multicolumn{1}{c|}{1.61}            & \multicolumn{1}{c|}{1.61}                               & \multicolumn{1}{c|}{3.26}             & \multicolumn{1}{c|}{3.64}                & 1.68           \\ \cline{2-12} 
	& 1          & \multicolumn{1}{c|}{1.18}            & \multicolumn{1}{c|}{1.40}               & \multicolumn{1}{c|}{0.25}             & \multicolumn{1}{c|}{0.93}                & 2.09           & \multicolumn{1}{c|}{0.36}            & \multicolumn{1}{c|}{0.36}                     & \multicolumn{1}{c|}{-0.92}            & \multicolumn{1}{c|}{-0.70}               & 1.70           \\ \cline{2-12} 
	& 1.6        & \multicolumn{1}{c|}{0.45}            & \multicolumn{1}{c|}{0.51}               &  \multicolumn{1}{c|}{-3.53}            & \multicolumn{1}{c|}{-3.13}               & 1.82           & \multicolumn{1}{c|}{0.01}            & \multicolumn{1}{c|}{0.01}                           & \multicolumn{1}{c|}{-2.28}            & \multicolumn{1}{c|}{-2.21}               & 0.87           \vspace{0.3cm}  \\ \hline
	\multirow{3}{*}{0.7} & 0.9        & \multicolumn{1}{c|}{2.30}            & \multicolumn{1}{c|}{2.70}               & \multicolumn{1}{c|}{5.25}             & \multicolumn{1}{c|}{6.18}                & 1.03           & \multicolumn{1}{c|}{1.98}            & \multicolumn{1}{c|}{2.20}                            & \multicolumn{1}{c|}{4.50}             & \multicolumn{1}{c|}{4.92}                & 1.21           \\ \cline{2-12} 
	& 1.4        & \multicolumn{1}{c|}{1.63}            & \multicolumn{1}{c|}{1.99}               &  \multicolumn{1}{c|}{2.13}             & \multicolumn{1}{c|}{3.25}                & 1.96           & \multicolumn{1}{c|}{0.63}            & \multicolumn{1}{c|}{0.77}                   & \multicolumn{1}{c|}{0.62}             & \multicolumn{1}{c|}{1.13}                & 1.91           \\ \cline{2-12} 
	& 2          & \multicolumn{1}{c|}{0.93}            & \multicolumn{1}{c|}{1.20}               & \multicolumn{1}{c|}{-0.90}            & \multicolumn{1}{c|}{0.02}                & 2.08           & \multicolumn{1}{c|}{0.17}            & \multicolumn{1}{c|}{0.23}                  & \multicolumn{1}{c|}{-1.49}            & \multicolumn{1}{c|}{-1.27}               & 1.48           \vspace{0.3cm} \\ \hline
	\multirow{3}{*}{1.2} & 1.5        & \multicolumn{1}{c|}{2.32}            & \multicolumn{1}{c|}{2.65}               &  \multicolumn{1}{c|}{5.38}             & \multicolumn{1}{c|}{6.17}                & 0.93           & \multicolumn{1}{c|}{2.07}            & \multicolumn{1}{c|}{2.31}                & \multicolumn{1}{c|}{4.68}             & \multicolumn{1}{c|}{5.11}                & 1.10           \\ \cline{2-12} 
	& 2          & \multicolumn{1}{c|}{1.97}            & \multicolumn{1}{c|}{2.36}               &  \multicolumn{1}{c|}{3.59}             & \multicolumn{1}{c|}{4.76}                & 1.69           & \multicolumn{1}{c|}{1.11}            & \multicolumn{1}{c|}{1.35}                            & \multicolumn{1}{c|}{2.23}             & \multicolumn{1}{c|}{2.95}                & 1.85           \\ \cline{2-12} 
	& 2.5        & \multicolumn{1}{c|}{1.55}            & \multicolumn{1}{c|}{1.93}               &  \multicolumn{1}{c|}{1.79}             & \multicolumn{1}{c|}{3.05}                & 1.20           & \multicolumn{1}{c|}{0.56}            & \multicolumn{1}{c|}{0.70}                 & \multicolumn{1}{c|}{0.28}             & \multicolumn{1}{c|}{0.85}                & 1.89         \vspace{0.3cm} \\ \hline\hline
  &        & \multicolumn{5}{c|}{$(n_1,n_2)=(6,8), ~~ (\mu_1, \mu_2)=(0.7,0.6)$}                                                                                                                                                                                                     & \multicolumn{5}{c}{$(n_1,n_2)=(9,4),  (\mu_1,\mu_2)=(0.7,0.6)$}  \vspace{0.3cm}                                                                                                                                                                                                    \\\hline\hline
	\multirow{3}{*}{0.4} & 0.6        & \multicolumn{1}{c|}{1.44}            & \multicolumn{1}{c|}{1.43}               & \multicolumn{1}{c|}{2.92}             & \multicolumn{1}{c|}{2.90}                & 1.80           & \multicolumn{1}{c|}{0.66}            & \multicolumn{1}{c|}{0.80}                  & \multicolumn{1}{c|}{0.65}             & \multicolumn{1}{c|}{0.88}                & 0.72           \\ \cline{2-12} 
	& 1          & \multicolumn{1}{c|}{0.31}            & \multicolumn{1}{c|}{0.31}               & \multicolumn{1}{c|}{-1.41}            & \multicolumn{1}{c|}{-1.41}               & 1.76           & \multicolumn{1}{c|}{0.47}            & \multicolumn{1}{c|}{0.42}                           & \multicolumn{1}{c|}{-0.34}            & \multicolumn{1}{c|}{-0.42}               & 0.77           \\ \cline{2-12} 
	& 1.6        & \multicolumn{1}{c|}{0.05}            & \multicolumn{1}{c|}{0.05}               & \multicolumn{1}{c|}{-3.02}            & \multicolumn{1}{c|}{-3.03}               & 0.98           & \multicolumn{1}{c|}{0.23}            & \multicolumn{1}{c|}{0.15}                    & \multicolumn{1}{c|}{-0.99}            & \multicolumn{1}{c|}{-1.18}               & 0.60           \vspace{0.3cm} \\ \hline
	\multirow{3}{*}{0.7} & 0.9        & \multicolumn{1}{c|}{1.96}            & \multicolumn{1}{c|}{2.20}               & \multicolumn{1}{c|}{3.86}             & \multicolumn{1}{c|}{4.31}                & 1.47           & \multicolumn{1}{c|}{0.71}            & \multicolumn{1}{c|}{1.05}                  & \multicolumn{1}{c|}{1.00}             & \multicolumn{1}{c|}{1.59}                & 0.63           \\ \cline{2-12} 
	& 1.4        & \multicolumn{1}{c|}{0.61}            & \multicolumn{1}{c|}{0.73}               & \multicolumn{1}{c|}{0.38}             & \multicolumn{1}{c|}{0.79}                & 1.96           & \multicolumn{1}{c|}{0.56}            & \multicolumn{1}{c|}{0.69}                 & \multicolumn{1}{c|}{0.09}             & \multicolumn{1}{c|}{0.39}                & 0.80           \\ \cline{2-12} 
	& 2          & \multicolumn{1}{c|}{0.19}            & \multicolumn{1}{c|}{0.23}               & \multicolumn{1}{c|}{-2.12}            & \multicolumn{1}{c|}{-1.89}               & 1.56           & \multicolumn{1}{c|}{0.42}            & \multicolumn{1}{c|}{0.44}                          & \multicolumn{1}{c|}{-0.55}            & \multicolumn{1}{c|}{-0.44}               & 0.74           \vspace{0.3cm} \\ \hline
	\multirow{3}{*}{1.2} & 1.5        & \multicolumn{1}{c|}{2.04}            & \multicolumn{1}{c|}{2.34}               & \multicolumn{1}{c|}{4.01}             & \multicolumn{1}{c|}{4.54}                & 1.38           & \multicolumn{1}{c|}{0.71}            & \multicolumn{1}{c|}{1.04}                     & \multicolumn{1}{c|}{1.05}             & \multicolumn{1}{c|}{1.67}                & 0.61           \\ \cline{2-12} 
	& 2          & \multicolumn{1}{c|}{1.07}            & \multicolumn{1}{c|}{1.38}               &  \multicolumn{1}{c|}{2.04}             & \multicolumn{1}{c|}{2.72}                & 1.93           & \multicolumn{1}{c|}{0.62}            & \multicolumn{1}{c|}{0.87}                         & \multicolumn{1}{c|}{0.44}             & \multicolumn{1}{c|}{0.95}                & 0.77           \\ \cline{2-12} 
	& 2.5        & \multicolumn{1}{c|}{0.53}            & \multicolumn{1}{c|}{0.70}        & \multicolumn{1}{c|}{0.01}             & \multicolumn{1}{c|}{0.67}                & 1.94           & \multicolumn{1}{c|}{0.55}            & \multicolumn{1}{c|}{0.71}                           & \multicolumn{1}{c|}{0.01}             & \multicolumn{1}{c|}{0.40}                & 0.80           \vspace{0.3cm} \\ \hline\hline
&        & \multicolumn{5}{c|}{$(n_1,n_2)=(10,15), (\mu_1,\mu_2)=(0.2, 0.7)$}                                                                                                                                                                                                      & \multicolumn{5}{c}{$(n_1,n_2)=(12,8), (\mu_1,\mu_2)=(0.2, 0.7)$}                                                                                                                                                                                                     \vspace{0.3cm} \\ \hline\hline
	\multirow{3}{*}{0.4} & 0.6        & \multicolumn{1}{c|}{1.10}            & \multicolumn{1}{c|}{1.03}               & \multicolumn{1}{c|}{1.51}             & \multicolumn{1}{c|}{1.43}                & 1.68           & \multicolumn{1}{c|}{0.65}            & \multicolumn{1}{c|}{0.58}               & \multicolumn{1}{c|}{1.17}             & \multicolumn{1}{c|}{1.05}                & 0.97           \\ \cline{2-12} 
	& 1          & \multicolumn{1}{c|}{0.09}            & \multicolumn{1}{c|}{0.08}               & \multicolumn{1}{c|}{-0.84}            & \multicolumn{1}{c|}{-0.86}               & 1.10           & \multicolumn{1}{c|}{0.13}            & \multicolumn{1}{c|}{0.09}                  & \multicolumn{1}{c|}{-0.35}            & \multicolumn{1}{c|}{-0.42}               & 0.76           \\ \cline{2-12} 
	& 1.6        & \multicolumn{1}{c|}{0.00}            & \multicolumn{1}{c|}{0.00}               & \multicolumn{1}{c|}{-1.28}            & \multicolumn{1}{c|}{-1.28}               & 0.33           & \multicolumn{1}{c|}{0.02}            & \multicolumn{1}{c|}{0.01}                 & \multicolumn{1}{c|}{-0.74}            & \multicolumn{1}{c|}{-0.77}               & 0.31           \vspace{0.3cm} \\ \hline
	\multirow{3}{*}{0.7} & 0.9        & \multicolumn{1}{c|}{1.54}            & \multicolumn{1}{c|}{1.64}               &  \multicolumn{1}{c|}{2.49}             & \multicolumn{1}{c|}{2.64}                & 1.45           & \multicolumn{1}{c|}{0.94}            & \multicolumn{1}{c|}{1.06}                      & \multicolumn{1}{c|}{1.66}             & \multicolumn{1}{c|}{1.86}                & 0.82           \\ \cline{2-12} 
	& 1.4        & \multicolumn{1}{c|}{0.37}            & \multicolumn{1}{c|}{0.40}               & \multicolumn{1}{c|}{-0.13}            & \multicolumn{1}{c|}{-0.04}               & 1.50           & \multicolumn{1}{c|}{0.28}            & \multicolumn{1}{c|}{0.31}                      & \multicolumn{1}{c|}{0.16}             & \multicolumn{1}{c|}{0.24}                & 0.95           \\ \cline{2-12} 
	& 2          & \multicolumn{1}{c|}{0.02}            & \multicolumn{1}{c|}{0.03}               &  \multicolumn{1}{c|}{-1.09}            & \multicolumn{1}{c|}{-1.06}               & 0.84           & \multicolumn{1}{c|}{0.08}            & \multicolumn{1}{c|}{0.07}                    & \multicolumn{1}{c|}{-0.52}            & \multicolumn{1}{c|}{-0.51}               & 0.62           \vspace{0.3cm} \\ \hline
	\multirow{3}{*}{1.2} & 1.5        & \multicolumn{1}{c|}{1.62}            & \multicolumn{1}{c|}{1.78}               & \multicolumn{1}{c|}{2.63}             & \multicolumn{1}{c|}{2.84}                & 1.38           & \multicolumn{1}{c|}{0.98}            & \multicolumn{1}{c|}{1.18}                       & \multicolumn{1}{c|}{1.74}             & \multicolumn{1}{c|}{2.03}                & 0.78           \\ \cline{2-12} 
	& 2          & \multicolumn{1}{c|}{0.75}            & \multicolumn{1}{c|}{0.88}               & \multicolumn{1}{c|}{0.89}             & \multicolumn{1}{c|}{1.13}                & 1.69           & \multicolumn{1}{c|}{0.49}            & \multicolumn{1}{c|}{0.62}                     & \multicolumn{1}{c|}{0.79}             & \multicolumn{1}{c|}{1.08}                & 1.00           \\ \cline{2-12} 
	& 2.5        & \multicolumn{1}{c|}{0.31}            & \multicolumn{1}{c|}{0.37}               & \multicolumn{1}{c|}{-0.30}            & \multicolumn{1}{c|}{-0.15}               & 1.44           & \multicolumn{1}{c|}{0.24}            & \multicolumn{1}{c|}{0.31}                     & \multicolumn{1}{c|}{0.04}             & \multicolumn{1}{c|}{0.23}                & 0.92           \\ \hline\hline
\end{tabular}
\end{table}
\begin{table}[h!]
	\centering 
	\caption{Percentage risk improvement of the improved estimators of $\mu_1$ over the BAEE under quadratic loss function $L_1(t)$}
	\vspace{0.2cm}
	\label{tab1_mu_1l4}
		\begin{tabular}{ccccccc|ccccc}
			\hline \hline\vspace{0.3cm}  
			$\sigma_1$           & $\sigma_2$ & \multicolumn{1}{c}{$\delta_{1s}^Q$} & \multicolumn{1}{c}{$\delta_{1S}^{*Q}$} & \multicolumn{1}{c}{$\delta_{1MV}^Q$} & \multicolumn{1}{c}{$\delta_{1MV}^{2Q}$} & $\delta_{1BZ}$ & \multicolumn{1}{c}{$\delta_{1s}^Q$} & \multicolumn{1}{c}{$\delta_{1S}^{*Q}$} &  \multicolumn{1}{c}{$\delta_{1MV}^Q$} & \multicolumn{1}{c}{$\delta_{1MV}^{2Q}$} & $\delta_{1BZ}$  \\\hline\hline \vspace{0.3cm}  
			&        &   \multicolumn{5}{c|}{$(n_1,n_2)=(8,15),   ~~(\mu_1,\mu_2)=(2,3)$}                                                                                                                                                                                                     & \multicolumn{5}{c}{$(n_1,n_2)=(10,17), (\mu_1,\mu_2)=(2,3)$}                                                                                                                                                                                                  \\ \hline\hline
		\multicolumn{1}{c|}{\multirow{3}{*}{0.4}} & \multicolumn{1}{c|}{0.6}        & \multicolumn{1}{c|}{1.41}            & \multicolumn{1}{c|}{0.51}               & \multicolumn{1}{c|}{2.36}               & \multicolumn{1}{c|}{0.11}                & 1.89           & \multicolumn{1}{c|}{0.84}            & \multicolumn{1}{c|}{0.17}               & \multicolumn{1}{c|}{1.37}               & \multicolumn{1}{c|}{-0.21}               & 1.55           \\ \cline{2-12} 
		\multicolumn{1}{c|}{}                     & \multicolumn{1}{c|}{1}          & \multicolumn{1}{c|}{0.24}            & \multicolumn{1}{c|}{0.08}               & \multicolumn{1}{c|}{-0.81}              & \multicolumn{1}{c|}{-1.38}               & 1.50           & \multicolumn{1}{c|}{0.00}            & \multicolumn{1}{c|}{0.02}              & \multicolumn{1}{c|}{-0.93}              & \multicolumn{1}{c|}{-1.10}               & 0.94           \\ \cline{2-12} 
		\multicolumn{1}{c|}{}                     & \multicolumn{1}{c|}{1.6}        & \multicolumn{1}{c|}{0.01}            & \multicolumn{1}{c|}{0.00}               & \multicolumn{1}{c|}{-1.62}              & \multicolumn{1}{c|}{-1.69}               & 0.59           & \multicolumn{1}{c|}{0.00}            & \multicolumn{1}{c|}{0.00}               & \multicolumn{1}{c|}{-1.13}              & \multicolumn{1}{c|}{-1.13}               & 0.23           \\ \hline
		\multicolumn{1}{c|}{\multirow{3}{*}{0.7}} & \multicolumn{1}{c|}{0.9}        & \multicolumn{1}{c|}{1.95}            & \multicolumn{1}{c|}{1.35}               & \multicolumn{1}{c|}{3.52}               & \multicolumn{1}{c|}{2.36}                & 1.51           & \multicolumn{1}{c|}{1.35}            & \multicolumn{1}{c|}{0.88}               & \multicolumn{1}{c|}{2.43}               & \multicolumn{1}{c|}{1.52}                & 1.30           \\ \cline{2-12} 
		\multicolumn{1}{c|}{}                     & \multicolumn{1}{c|}{1.4}        & \multicolumn{1}{c|}{0.61}            & \multicolumn{1}{c|}{0.35}               & \multicolumn{1}{c|}{0.19}               & \multicolumn{1}{c|}{-0.42}               & 1.89           & \multicolumn{1}{c|}{0.17}            & \multicolumn{1}{c|}{0.06}               & \multicolumn{1}{c|}{-0.30}              & \multicolumn{1}{c|}{-0.68}               & 1.36           \\ \cline{2-12} 
		\multicolumn{1}{c|}{}                     & \multicolumn{1}{c|}{2}          & \multicolumn{1}{c|}{0.13}            & \multicolumn{1}{c|}{0.07}               & \multicolumn{1}{c|}{-1.21}              & \multicolumn{1}{c|}{-1.38}               & 1.22           & \multicolumn{1}{c|}{0.02}           & \multicolumn{1}{c|}{0.02}              & \multicolumn{1}{c|}{-1.07}              & \multicolumn{1}{c|}{-1.12}               & 0.69           \\ \hline
		\multicolumn{1}{c|}{\multirow{3}{*}{1.2}} & \multicolumn{1}{c|}{1.5}        & \multicolumn{1}{c|}{2.04}            & \multicolumn{1}{c|}{1.84}               & \multicolumn{1}{c|}{3.69}               & \multicolumn{1}{c|}{3.35}                & 1.41           & \multicolumn{1}{c|}{1.45}            & \multicolumn{1}{c|}{1.27}               & \multicolumn{1}{c|}{2.60}               & \multicolumn{1}{c|}{2.30}                & 1.22           \\ \cline{2-12} 
		\multicolumn{1}{c|}{}                     & \multicolumn{1}{c|}{2}          & \multicolumn{1}{c|}{1.05}            & \multicolumn{1}{c|}{0.89}               & \multicolumn{1}{c|}{1.53}               & \multicolumn{1}{c|}{1.13}                & 1.98           & \multicolumn{1}{c|}{0.53}            & \multicolumn{1}{c|}{0.41}               & \multicolumn{1}{c|}{0.67}               & \multicolumn{1}{c|}{0.40}                & 1.55           \\ \cline{2-12} 
		\multicolumn{1}{c|}{}                     & \multicolumn{1}{c|}{2.5}        & \multicolumn{1}{c|}{0.52}            & \multicolumn{1}{c|}{0.42}               & \multicolumn{1}{c|}{-0.03}              & \multicolumn{1}{c|}{-0.25}               & 1.83           & \multicolumn{1}{c|}{0.12}            & \multicolumn{1}{c|}{0.08}               & \multicolumn{1}{c|}{-0.45}              & \multicolumn{1}{c|}{-0.58}               & 1.29    
		\vspace{0.3cm} \\ \hline\hline
		&        & \multicolumn{5}{c|}{$(n_1,n_2)=(14,20), ~~ (\mu_1, \mu_2)=(3,5)$}                                                                                                                                                                                                     & \multicolumn{5}{c}{$(n_1,n_2)=(20,25),  (\mu_1,\mu_2)=(3,5)$}  \vspace{0.3cm}                                                                                                                                                                                                    \\\hline\hline
		
		\multicolumn{1}{c|}{\multirow{3}{*}{0.4}} & \multicolumn{1}{c|}{0.6}        & \multicolumn{1}{c|}{0.58}            & \multicolumn{1}{c|}{0.03}               & \multicolumn{1}{c|}{0.99}               & \multicolumn{1}{c|}{-0.29}               & 1.10           & \multicolumn{1}{c|}{0.36}            & \multicolumn{1}{c|}{0.02}               & \multicolumn{1}{c|}{0.36}               & \multicolumn{1}{c|}{-0.32}               & 0.95           \\ \cline{2-12} 
		\multicolumn{1}{c|}{}                     & \multicolumn{1}{c|}{1}          & \multicolumn{1}{c|}{0.02}            & \multicolumn{1}{c|}{0.00}               & \multicolumn{1}{c|}{-0.35}              & \multicolumn{1}{c|}{-0.43}               & 0.53           & \multicolumn{1}{c|}{0.01}            & \multicolumn{1}{c|}{0.00}               & \multicolumn{1}{c|}{-0.36}              & \multicolumn{1}{c|}{-0.38}               & 0.26           \\ \cline{2-12} 
		\multicolumn{1}{c|}{}                     & \multicolumn{1}{c|}{1.6}        & \multicolumn{1}{c|}{0.00}            & \multicolumn{1}{c|}{0.00}               & \multicolumn{1}{c|}{-0.44}              & \multicolumn{1}{c|}{-0.44}               & 0.09           & \multicolumn{1}{c|}{0.00}            & \multicolumn{1}{c|}{0.00}               & \multicolumn{1}{c|}{-0.38}              & \multicolumn{1}{c|}{-0.38}               & 0.02           \\ \hline
		\multicolumn{1}{c|}{\multirow{3}{*}{0.7}} & \multicolumn{1}{c|}{0.9}        & \multicolumn{1}{c|}{1.03}            & \multicolumn{1}{c|}{0.36}               & \multicolumn{1}{c|}{1.70}               & \multicolumn{1}{c|}{0.59}                & 0.91           & \multicolumn{1}{c|}{0.70}            & \multicolumn{1}{c|}{0.26}               & \multicolumn{1}{c|}{0.95}               & \multicolumn{1}{c|}{0.22}                & 0.95           \\ \cline{2-12} 
		\multicolumn{1}{c|}{}                     & \multicolumn{1}{c|}{1.4}        & \multicolumn{1}{c|}{0.11}            & \multicolumn{1}{c|}{0.02}               & \multicolumn{1}{c|}{-0.08}              & \multicolumn{1}{c|}{-0.33}               & 0.88           & \multicolumn{1}{c|}{0.04}            & \multicolumn{1}{c|}{0.01}               & \multicolumn{1}{c|}{-0.27}              & \multicolumn{1}{c|}{-0.35}               & 0.56           \\ \cline{2-12} 
		\multicolumn{1}{c|}{}                     & \multicolumn{1}{c|}{2}          & \multicolumn{1}{c|}{0.00}            & \multicolumn{1}{c|}{0.00}               & \multicolumn{1}{c|}{-0.41}              & \multicolumn{1}{c|}{-0.43}               & 0.35           & \multicolumn{1}{c|}{0.00}            & \multicolumn{1}{c|}{0.00}               & \multicolumn{1}{c|}{-0.37}              & \multicolumn{1}{c|}{-0.38}               & 0.14           \\ \hline
		\multicolumn{1}{c|}{\multirow{3}{*}{1.2}} & \multicolumn{1}{c|}{1.5}        & \multicolumn{1}{c|}{1.10}            & \multicolumn{1}{c|}{0.75}               & \multicolumn{1}{c|}{1.82}               & \multicolumn{1}{c|}{1.31}                & 0.85           & \multicolumn{1}{c|}{0.77}            & \multicolumn{1}{c|}{0.53}               & \multicolumn{1}{c|}{1.06}               & \multicolumn{1}{c|}{0.68}                & 0.93           \\ \cline{2-12} 
		\multicolumn{1}{c|}{}                     & \multicolumn{1}{c|}{2}          & \multicolumn{1}{c|}{0.34}            & \multicolumn{1}{c|}{0.19}               & \multicolumn{1}{c|}{0.50}               & \multicolumn{1}{c|}{0.12}                & 1.08           & \multicolumn{1}{c|}{0.18}            & \multicolumn{1}{c|}{0.08}               & \multicolumn{1}{c|}{0.06}               & \multicolumn{1}{c|}{-0.13}               & 0.84           \\ \cline{2-12} 
		\multicolumn{1}{c|}{}                     & \multicolumn{1}{c|}{2.5}        & \multicolumn{1}{c|}{0.09}            & \multicolumn{1}{c|}{0.04}               & \multicolumn{1}{c|}{-0.15}              & \multicolumn{1}{c|}{-0.28}               & 0.81           & \multicolumn{1}{c|}{0.03}            & \multicolumn{1}{c|}{0.02}               & \multicolumn{1}{c|}{-0.30}              & \multicolumn{1}{c|}{-0.34}               & 0.50          
		\vspace{0.3cm} \\ \hline\hline
		&        & \multicolumn{5}{c|}{$(n_1,n_2)=(27,35), ~~ (\mu_1, \mu_2)=(5,7)$}                                                                                                                                                                                                     & \multicolumn{5}{c}{$(n_1,n_2)=(32,38),  (\mu_1,\mu_2)=(5,7)$}  \vspace{0.3cm}                                                                                                                                                                                                    \\\hline\hline

		\multicolumn{1}{c|}{\multirow{3}{*}{0.4}} & \multicolumn{1}{c|}{0.6}        & \multicolumn{1}{c|}{0.15}            & \multicolumn{1}{c|}{0.00}               & \multicolumn{1}{c|}{0.20}               & \multicolumn{1}{c|}{-0.09}               & 0.53           & \multicolumn{1}{c|}{0.12}            & \multicolumn{1}{c|}{0.00}               & \multicolumn{1}{c|}{0.16}               & \multicolumn{1}{c|}{-0.06}               & 0.48           \\ \cline{2-12} 
		\multicolumn{1}{c|}{}                     & \multicolumn{1}{c|}{1}          & \multicolumn{1}{c|}{0.00}            & \multicolumn{1}{c|}{0.00}               & \multicolumn{1}{c|}{-0.09}              & \multicolumn{1}{c|}{-0.10}               & 0.08           & \multicolumn{1}{c|}{0.00}            & \multicolumn{1}{c|}{0.00}               & \multicolumn{1}{c|}{-0.07}              & \multicolumn{1}{c|}{-0.07}               & 0.05           \\ \cline{2-12} 
		\multicolumn{1}{c|}{}                     & \multicolumn{1}{c|}{1.6}        & \multicolumn{1}{c|}{0.00}            & \multicolumn{1}{c|}{0.00}               & \multicolumn{1}{c|}{-0.10}              & \multicolumn{1}{c|}{-0.10}               & 0.00           & \multicolumn{1}{c|}{0.00}            & \multicolumn{1}{c|}{0.00}               & \multicolumn{1}{c|}{-0.07}              & \multicolumn{1}{c|}{-0.07}               & 0.00           \\ \hline
		\multicolumn{1}{c|}{\multirow{3}{*}{0.7}} & \multicolumn{1}{c|}{0.9}        & \multicolumn{1}{c|}{0.35}            & \multicolumn{1}{c|}{0.10}               & \multicolumn{1}{c|}{0.55}               & \multicolumn{1}{c|}{0.12}                & 0.53           & \multicolumn{1}{c|}{0.35}            & \multicolumn{1}{c|}{0.07}               & \multicolumn{1}{c|}{0.51}               & \multicolumn{1}{c|}{0.09}                & 0.52           \\ \cline{2-12} 
		\multicolumn{1}{c|}{}                     & \multicolumn{1}{c|}{1.4}        & \multicolumn{1}{c|}{0.01}            & \multicolumn{1}{c|}{0.00}               & \multicolumn{1}{c|}{-0.08}              & \multicolumn{1}{c|}{-0.09}               & 0.24           & \multicolumn{1}{c|}{0.00}           & \multicolumn{1}{c|}{0.00}               & \multicolumn{1}{c|}{-0.06}              & \multicolumn{1}{c|}{-0.07}               & 0.19           \\ \cline{2-12} 
		\multicolumn{1}{c|}{}                     & \multicolumn{1}{c|}{2}          & \multicolumn{1}{c|}{0.00}            & \multicolumn{1}{c|}{0.00}               & \multicolumn{1}{c|}{-0.10}              & \multicolumn{1}{c|}{-0.10}               & 0.03           & \multicolumn{1}{c|}{0.00}            & \multicolumn{1}{c|}{0.00}               & \multicolumn{1}{c|}{-0.07}              & \multicolumn{1}{c|}{-0.07}               & 0.02           \\ \hline
		\multicolumn{1}{c|}{\multirow{3}{*}{1.2}} & \multicolumn{1}{c|}{1.5}        & \multicolumn{1}{c|}{0.40}            & \multicolumn{1}{c|}{0.24}               & \multicolumn{1}{c|}{0.64}               & \multicolumn{1}{c|}{0.34}                & 0.51           & \multicolumn{1}{c|}{0.40}            & \multicolumn{1}{c|}{0.21}               & \multicolumn{1}{c|}{0.59}               & \multicolumn{1}{c|}{0.32}                & 0.51           \\ \cline{2-12} 
		\multicolumn{1}{c|}{}                     & \multicolumn{1}{c|}{2}          & \multicolumn{1}{c|}{0.06}            & \multicolumn{1}{c|}{0.01}               & \multicolumn{1}{c|}{0.05}               & \multicolumn{1}{c|}{-0.05}               & 0.44           & \multicolumn{1}{c|}{0.04}            & \multicolumn{1}{c|}{0.01}               & \multicolumn{1}{c|}{0.02}               & \multicolumn{1}{c|}{-0.04}               & 0.38           \\ \cline{2-12} 
		\multicolumn{1}{c|}{}                     & \multicolumn{1}{c|}{2.5}        & \multicolumn{1}{c|}{0.01}            & \multicolumn{1}{c|}{0.00}               & \multicolumn{1}{c|}{-0.09}              & \multicolumn{1}{c|}{-0.09}               & 0.20           & \multicolumn{1}{c|}{0.00}           & \multicolumn{1}{c|}{0.00}              & \multicolumn{1}{c|}{-0.07}              & \multicolumn{1}{c|}{-0.07}               & 0.15           \\ \hline \hline
	\end{tabular}
\end{table}
\begin{table}[h!]
		\centering
	\caption{Percentage risk improvement of $\delta_{1MV}^Q$ and $\delta_{1MV}^{2Q}$ of $\mu_1$ over the UMVUE under quadratic loss function $L_1(t)$}
	\vspace{0.2cm}
	\label{tab1_mu_1l2}
	\begin{tabular}{|c|cccc|cc|}
		\hline
		($\mu_1,\mu_2$)             & \multicolumn{1}{c|}{$\sigma_1$}           & \multicolumn{1}{c|}{$\sigma_2$} & \multicolumn{1}{c|}{$\delta_{1MV}^Q$} & $\delta_{1MV}^{2Q}$ & \multicolumn{1}{c|}{$\delta_{1MV}^Q$} & $\delta_{1MV}^{2Q}$ \\ \hline
		\multirow{10}{*}{(0.1,0.3)} & \multicolumn{4}{c|}{$(n_1,n_2)=(6,11)$}                                                                                                   & \multicolumn{2}{c|}{$ (n_1,n_2)=(4,5)$}                     \\ \cline{2-7} 
		& \multicolumn{1}{c|}{\multirow{3}{*}{0.4}} & \multicolumn{1}{c|}{0.6}        & \multicolumn{1}{c|}{5.965036}         & 6.351358            & \multicolumn{1}{c|}{9.94595}          & 10.90282            \\ \cline{3-7} 
		& \multicolumn{1}{c|}{}                     & \multicolumn{1}{c|}{1}          & \multicolumn{1}{c|}{2.005833}         & 2.265678            & \multicolumn{1}{c|}{5.52922}          & 6.379911            \\ \cline{3-7} 
		& \multicolumn{1}{c|}{}                     & \multicolumn{1}{c|}{1.6}        & \multicolumn{1}{c|}{0.3267855}        & 0.4265232           & \multicolumn{1}{c|}{2.369576}         & 2.688524            \\ \cline{2-7} 
		& \multicolumn{1}{c|}{\multirow{3}{*}{0.7}} & \multicolumn{1}{c|}{0.9}        & \multicolumn{1}{c|}{7.087558}         & 7.535222            & \multicolumn{1}{c|}{10.85123}         & 11.84486            \\ \cline{3-7} 
		& \multicolumn{1}{c|}{}                     & \multicolumn{1}{c|}{1.4}        & \multicolumn{1}{c|}{3.542498}         & 4.041664            & \multicolumn{1}{c|}{7.61183}          & 8.811912            \\ \cline{3-7} 
		& \multicolumn{1}{c|}{}                     & \multicolumn{1}{c|}{2}          & \multicolumn{1}{c|}{1.412706}         & 1.670769            & \multicolumn{1}{c|}{4.44379}          & 5.36237             \\ \cline{2-7} 
		& \multicolumn{1}{c|}{\multirow{3}{*}{1.2}} & \multicolumn{1}{c|}{1.5}        & \multicolumn{1}{c|}{7.262581}         & 7.714264            & \multicolumn{1}{c|}{10.98141}         & 11.85648            \\ \cline{3-7} 
		& \multicolumn{1}{c|}{}                     & \multicolumn{1}{c|}{2}          & \multicolumn{1}{c|}{5.061812}         & 5.709763            & \multicolumn{1}{c|}{9.172332}         & 10.44139            \\ \cline{3-7} 
		& \multicolumn{1}{c|}{}                     & \multicolumn{1}{c|}{2.5}        & \multicolumn{1}{c|}{3.224353}         & 3.779561            & \multicolumn{1}{c|}{7.239574}         & 8.597475            \\ \hline
		\multirow{10}{*}{(0.7,0.6)} & \multicolumn{4}{c|}{$(n_1,n_2)=(6,8)$}                                                                                                    & \multicolumn{2}{c|}{$(n_1,n_2)=(9,4)$}                            \\ \cline{2-7} 
		& \multicolumn{1}{c|}{\multirow{3}{*}{0.4}} & \multicolumn{1}{c|}{0.6}        & \multicolumn{1}{c|}{5.897498}         & 5.959126            & \multicolumn{1}{c|}{2.597914}         & 2.836676            \\ \cline{3-7} 
		& \multicolumn{1}{c|}{}                     & \multicolumn{1}{c|}{1}          & \multicolumn{1}{c|}{2.421264}         & 2.435001            & \multicolumn{1}{c|}{1.510962}         & 1.432583            \\ \cline{3-7} 
		& \multicolumn{1}{c|}{}                     & \multicolumn{1}{c|}{1.6}        & \multicolumn{1}{c|}{0.6485672}        & 0.672182            & \multicolumn{1}{c|}{0.779326}         & 0.668592            \\ \cline{2-7} 
		& \multicolumn{1}{c|}{\multirow{3}{*}{0.7}} & \multicolumn{1}{c|}{0.9}        & \multicolumn{1}{c|}{6.850938}         & 7.274175            & \multicolumn{1}{c|}{2.82484}          & 3.46123             \\ \cline{3-7} 
		& \multicolumn{1}{c|}{}                     & \multicolumn{1}{c|}{1.4}        & \multicolumn{1}{c|}{3.831223}         & 4.209723            & \multicolumn{1}{c|}{1.983337}         & 2.313678            \\ \cline{3-7} 
		& \multicolumn{1}{c|}{}                     & \multicolumn{1}{c|}{2}          & \multicolumn{1}{c|}{1.796149}         & 1.996802            & \multicolumn{1}{c|}{1.240297}         & 1.365613            \\ \cline{2-7} 
		& \multicolumn{1}{c|}{\multirow{3}{*}{1.2}} & \multicolumn{1}{c|}{1.5}        & \multicolumn{1}{c|}{6.989729}         & 7.502268            & \multicolumn{1}{c|}{6.989729}         & 7.502268            \\ \cline{3-7} 
		& \multicolumn{1}{c|}{}                     & \multicolumn{1}{c|}{2}          & \multicolumn{1}{c|}{5.156025}         & 5.771038            & \multicolumn{1}{c|}{5.156025}         & 5.771038            \\ \cline{3-7} 
		& \multicolumn{1}{c|}{}                     & \multicolumn{1}{c|}{2.5}        & \multicolumn{1}{c|}{3.553044}         & 4.105228            & \multicolumn{1}{c|}{3.553044}         & 4.105228            \\ \hline
		\multirow{10}{*}{(0.2,0.7)} & \multicolumn{4}{c|}{$(n_1,n_2)=(10,15)$}                                                                                                  & \multicolumn{2}{c|}{$(n_1,n_2)=(12,8)$}                     \\ \cline{2-7} 
		& \multicolumn{1}{c|}{0.4}                  & \multicolumn{1}{c|}{0.6}        & \multicolumn{1}{c|}{2.664594}         & 2.569208            & \multicolumn{1}{c|}{1.797314}         & 1.733149            \\ \cline{2-7} 
		& \multicolumn{1}{c|}{}                     & \multicolumn{1}{c|}{1}          & \multicolumn{1}{c|}{0.5152419}        & 0.4904035           & \multicolumn{1}{c|}{0.562659}         & 0.509239            \\ \cline{2-7} 
		& \multicolumn{1}{c|}{}                     & \multicolumn{1}{c|}{1.6}        & \multicolumn{1}{c|}{0.02423753}       & 0.02200             & \multicolumn{1}{c|}{0.104269}         & 0.08634             \\ \cline{2-7} 
		& \multicolumn{1}{c|}{\multirow{3}{*}{0.7}} & \multicolumn{1}{c|}{0.9}        & \multicolumn{1}{c|}{3.4495}           & 3.575824            & \multicolumn{1}{c|}{2.172759}         & 2.393397            \\ \cline{3-7} 
		& \multicolumn{1}{c|}{}                     & \multicolumn{1}{c|}{1.4}        & \multicolumn{1}{c|}{1.194514}         & 1.282934            & \multicolumn{1}{c|}{1.001425}         & 1.103839            \\ \cline{3-7} 
		& \multicolumn{1}{c|}{}                     & \multicolumn{1}{c|}{2}          & \multicolumn{1}{c|}{0.2323653}        & 0.2696534           & \multicolumn{1}{c|}{0.366258}         & 0.406029            \\ \cline{2-7} 
		& \multicolumn{1}{c|}{\multirow{3}{*}{1.2}} & \multicolumn{1}{c|}{1.5}        & \multicolumn{1}{c|}{3.576267}         & 3.78095             & \multicolumn{1}{c|}{2.236442}         & 2.498017            \\ \cline{3-7} 
		& \multicolumn{1}{c|}{}                     & \multicolumn{1}{c|}{2}          & \multicolumn{1}{c|}{2.074687}         & 2.291207            & \multicolumn{1}{c|}{1.48783}          & 1.722845            \\ \cline{3-7} 
		& \multicolumn{1}{c|}{}                     & \multicolumn{1}{c|}{2.5}        & \multicolumn{1}{c|}{1.04078}          & 1.18672             & \multicolumn{1}{c|}{0.90751}          & 1.086625            \\ \hline
	\end{tabular}
\end{table}

\begin{table}[h!]
	\centering
	\caption{Percentage risk improvement of $\delta_{1MV}^Q$ and $\delta_{1MV}^{2Q}$ of $\mu_1$ over the UMVUE under quadratic loss function $L_1(t)$}
	\vspace{0.2cm}
	\label{tab1_mu_1l3}
		\begin{tabular}{|c|cccc|cc|}
			\hline
			$(\mu_1,\mu_2)$         & \multicolumn{1}{c|}{$\sigma_1$}           & \multicolumn{1}{c|}{$\sigma_2$} & \multicolumn{1}{c|}{$\delta_{1MV}^{Q}$} & $\delta_{1MV}^{2Q}$ & \multicolumn{1}{c|}{$\delta_{1MV}^{Q}$} & $\delta_{1MV}^{2Q}$ \\ \hline
			\multirow{10}{*}{(2,3)} & \multicolumn{4}{c|}{$(n_1,n_2)=(8,15)$}                                                                                                     & \multicolumn{2}{c|}{$(n_1,n_2)=(10,17)$}                      \\ \cline{2-7} 
			& \multicolumn{1}{c|}{\multirow{3}{*}{0.4}} & \multicolumn{1}{c|}{0.6}        & \multicolumn{1}{c|}{4.100092}           & 1.982114            & \multicolumn{1}{c|}{2.653808}           & 1.209291            \\ \cline{3-7} 
			& \multicolumn{1}{c|}{}                     & \multicolumn{1}{c|}{1}          & \multicolumn{1}{c|}{0.987338}           & 0.321735            & \multicolumn{1}{c|}{0.444077}           & 0.117267            \\ \cline{3-7} 
			& \multicolumn{1}{c|}{}                     & \multicolumn{1}{c|}{1.6}        & \multicolumn{1}{c|}{0.068328}           & 0.010054            & \multicolumn{1}{c|}{0.013806}           & 0                   \\ \cline{2-7} 
			& \multicolumn{1}{c|}{\multirow{3}{*}{0.7}} & \multicolumn{1}{c|}{0.9}        & \multicolumn{1}{c|}{5.15343}            & 4.126812            & \multicolumn{1}{c|}{3.543086}           & 2.765168            \\ \cline{3-7} 
			& \multicolumn{1}{c|}{}                     & \multicolumn{1}{c|}{1.4}        & \multicolumn{1}{c|}{2.05438}            & 1.402398            & \multicolumn{1}{c|}{1.117098}           & 0.754051            \\ \cline{3-7} 
			& \multicolumn{1}{c|}{}                     & \multicolumn{1}{c|}{2}          & \multicolumn{1}{c|}{0.530457}           & 0.307990            & \multicolumn{1}{c|}{0.182608}           & 0.106249            \\ \cline{2-7} 
			& \multicolumn{1}{c|}{\multirow{3}{*}{1.2}} & \multicolumn{1}{c|}{1.5}        & \multicolumn{1}{c|}{5.320618}           & 5.015111            & \multicolumn{1}{c|}{3.676835}           & 3.454983            \\ \cline{3-7} 
			& \multicolumn{1}{c|}{}                     & \multicolumn{1}{c|}{2}          & \multicolumn{1}{c|}{3.298495}           & 2.941397            & \multicolumn{1}{c|}{2.028779}           & 1.778652            \\ \cline{3-7} 
			& \multicolumn{1}{c|}{}                     & \multicolumn{1}{c|}{2.5}        & \multicolumn{1}{c|}{1.81217}            & 1.589339            & \multicolumn{1}{c|}{0.968497}           & 0.8453503           \\ \hline
			& \multicolumn{4}{c|}{$(n_1,n_2)=(14,20)$}                                                                                                    & \multicolumn{2}{c|}{$(n_1,n_2)=(20,25)$}                      \\ \hline
			\multirow{9}{*}{(3,5)}  & \multicolumn{1}{c|}{\multirow{3}{*}{0.4}} & \multicolumn{1}{c|}{0.6}        & \multicolumn{1}{c|}{1.495304}           & 0.2005839           & \multicolumn{1}{c|}{0.815367}           & 0.095532            \\ \cline{3-7} 
			& \multicolumn{1}{c|}{}                     & \multicolumn{1}{c|}{1}          & \multicolumn{1}{c|}{0.114227}           & 0                   & \multicolumn{1}{c|}{0.025428}           & 0                   \\ \cline{3-7} 
			& \multicolumn{1}{c|}{}                     & \multicolumn{1}{c|}{1.6}        & \multicolumn{1}{c|}{0}                  & 0                   & \multicolumn{1}{c|}{0}                  & 0                   \\ \cline{2-7} 
			& \multicolumn{1}{c|}{\multirow{3}{*}{0.7}} & \multicolumn{1}{c|}{0.9}        & \multicolumn{1}{c|}{2.214141}           & 1.128142            & \multicolumn{1}{c|}{1.404343}           & 0.681709            \\ \cline{3-7} 
			& \multicolumn{1}{c|}{}                     & \multicolumn{1}{c|}{1.4}        & \multicolumn{1}{c|}{0.479748}           & 0.132767            & \multicolumn{1}{c|}{0.172804}           & 0.044693            \\ \cline{3-7} 
			& \multicolumn{1}{c|}{}                     & \multicolumn{1}{c|}{2}          & \multicolumn{1}{c|}{0.042714}           & 0.002477            & \multicolumn{1}{c|}{0.001427}           & 0                   \\ \cline{2-7} 
			& \multicolumn{1}{c|}{\multirow{3}{*}{1.2}} & \multicolumn{1}{c|}{1.5}        & \multicolumn{1}{c|}{2.332249}           & 1.824857            & \multicolumn{1}{c|}{1.506951}           & 1.147355            \\ \cline{3-7} 
			& \multicolumn{1}{c|}{}                     & \multicolumn{1}{c|}{2}          & \multicolumn{1}{c|}{1.038007}           & 0.678675            & \multicolumn{1}{c|}{0.486567}           & 0.3184061           \\ \cline{3-7} 
			& \multicolumn{1}{c|}{}                     & \multicolumn{1}{c|}{2.5}        & \multicolumn{1}{c|}{0.386580}           & 0.219414            & \multicolumn{1}{c|}{0.124521}           & 0.065179            \\ \hline
			& \multicolumn{4}{c|}{$(n_1,n_2)=(27,35)$}                                                                                                    & \multicolumn{2}{c|}{$(n_1,n_2)=(32,38)$}                      \\ \hline
			\multirow{9}{*}{(5,7)}  & \multicolumn{1}{c|}{\multirow{3}{*}{0.4}} & \multicolumn{1}{c|}{0.6}        & \multicolumn{1}{c|}{0.377370}           & 0.013292            & \multicolumn{1}{c|}{0.262095}           & 0.008266            \\ \cline{3-7} 
			& \multicolumn{1}{c|}{}                     & \multicolumn{1}{c|}{1}          & \multicolumn{1}{c|}{0}                  & 0                   & \multicolumn{1}{c|}{0}                  & 0                   \\ \cline{3-7} 
			& \multicolumn{1}{c|}{}                     & \multicolumn{1}{c|}{1.6}        & \multicolumn{1}{c|}{0}                  & 0                   & \multicolumn{1}{c|}{0}                  & 0                   \\ \cline{2-7} 
			& \multicolumn{1}{c|}{\multirow{3}{*}{0.7}} & \multicolumn{1}{c|}{0.9}        & \multicolumn{1}{c|}{0.839457}           & 0.265718            & \multicolumn{1}{c|}{0.648369}           & 0.197522            \\ \cline{3-7} 
			& \multicolumn{1}{c|}{}                     & \multicolumn{1}{c|}{1.4}        & \multicolumn{1}{c|}{0.035650}           & 0.000806            & \multicolumn{1}{c|}{0.019406}           & 0                   \\ \cline{3-7} 
			& \multicolumn{1}{c|}{}                     & \multicolumn{1}{c|}{2}          & \multicolumn{1}{c|}{0}                  & 0                   & \multicolumn{1}{c|}{0}                  & 0                   \\ \cline{2-7} 
			& \multicolumn{1}{c|}{\multirow{3}{*}{1.2}} & \multicolumn{1}{c|}{1.5}        & \multicolumn{1}{c|}{0.933795}           & 0.578254            & \multicolumn{1}{c|}{0.731979}           & 0.441341            \\ \cline{3-7} 
			& \multicolumn{1}{c|}{}                     & \multicolumn{1}{c|}{2}          & \multicolumn{1}{c|}{0.196293}           & 0.083866            & \multicolumn{1}{c|}{0.128046}           & 0.047961            \\ \cline{3-7} 
			& \multicolumn{1}{c|}{}                     & \multicolumn{1}{c|}{2.5}        & \multicolumn{1}{c|}{0.024242}           & 0.005683            & \multicolumn{1}{c|}{0.010899}           & 0.000624            \\ \hline
		\end{tabular}
\end{table}

\clearpage
\vspace{5cm}

\section{Conclusion}\label{sec8}
In this manuscript, we have investigated component-wise estimation location parameters of two exponential distributions with ordered scale parameters under a bowl-shaped invariant loss function. 
The inadmissibility of the BAEE has been proved by proposing improved estimators. We have obtained two dominating \cite{stein1964} type estimators. Further, a class of dominating estimators has been presented using the IERD approach of \cite{kubokawa1994unified}. It is noted that the boundary estimator of this class is the Brewster-Zidek type improved estimator. We have proved that the Brewster-Zidek-type improved estimator is generalized Bayes. We have obtained the expressions of improved estimators for squared error loss and linex loss. We also consider the estimation problem with respect to the generalized Pitman closeness criterion. Under the generalized Pitman closeness criterion, a Stein-type improved estimator is proposed. The estimators enhanced through the pitman closeness criterion demonstrate robustness as they do not depend on the form of the loss function. Finally, we performed a simulation study to verify the paper's findings and assess the performance of different competing estimators.
\bibliography{bib_location}
\end{document}